\documentclass[a4paper]{amsart}
\usepackage[margin=1.25in]{geometry}
\usepackage[justification=centering]{caption}
\usepackage{float}
\usepackage{listings}
\lstdefinelanguage{Sage}[]{Python}
{morekeywords={False,sage,True},sensitive=true}
\usepackage[colorinlistoftodos]{todonotes}
\usepackage{amscd,amsmath,amsxtra,amsthm,amssymb,stmaryrd,xr,mathrsfs,mathtools,enumerate,commath, comment}
\usepackage{orcidlink}
\usepackage{stmaryrd}
\usepackage{xcolor}
\usepackage{commath}
\usepackage{comment}
\usepackage{tikz-cd}
\usepackage{enumitem}
\usepackage{longtable} 
\usepackage{pdflscape} 
\usepackage{booktabs}
\usepackage{hyperref}
\definecolor{dblackcolor}{rgb}{0.0,0.0,0.0}
\definecolor{dbluecolor}{rgb}{0.01,0.02,0.7}
\definecolor{dgreencolor}{rgb}{0.2,0.4,0.0}
\definecolor{dgraycolor}{rgb}{0.30,0.3,0.30}

\definecolor{backcolour}{rgb}{0.95,0.95,0.92}
\definecolor{vegasgold}{rgb}{0.77, 0.7, 0.35}
\definecolor{darkgoldenrod}{rgb}{0.72, 0.53, 0.04}
\definecolor{gold(metallic)}{rgb}{0.83, 0.69, 0.22}
\hypersetup{
 colorlinks=true,
 linkcolor=darkgoldenrod,
 filecolor=brown,      
 urlcolor=gold(metallic),
 citecolor=darkgoldenrod,
 pdftitle={On the number of subrings of $\mathbb{Z}^n$ of prime power index},
 }

\usepackage[all,cmtip]{xy}

\DeclareFontFamily{U}{wncy}{}
\DeclareFontShape{U}{wncy}{m}{n}{<->wncyr10}{}
\DeclareSymbolFont{mcy}{U}{wncy}{m}{n}
\DeclareMathSymbol{\Sh}{\mathord}{mcy}{"58}
\usepackage[T2A,T1]{fontenc}
\usepackage[OT2,T1]{fontenc}

\newtheorem{theorem}{Theorem}[section]
\newtheorem{lemma}[theorem]{Lemma}
\newtheorem{question}[theorem]{Question}
\newtheorem{conjecture}[theorem]{Conjecture}
\newtheorem*{theorem*}{Theorem}
\newtheorem*{corollary*}{Corollary}

\newtheorem{proposition}[theorem]{Proposition}

\newtheorem{definition}[theorem]{Definition}

\numberwithin{equation}{section}

\theoremstyle{remark}
\newtheorem{remark}[theorem]{Remark}

\newcommand{\Z}{\mathbb{Z}}

\newcommand{\Q}{\mathbb{Q}}
\newcommand{\F}{\mathbb{F}}

\newcommand{\cO}{\mathcal{O}}

\renewcommand{\Col}{\mathrm{Col}}

\newcommand{\op}[1]{\operatorname{#1}}

\begin{document}
\title[Subrings of $\mathbb{Z}^n$ of prime power index]{On the number of subrings of $\Z^n$ of prime power index}

\author[H.~Mishra]{Hrishabh Mishra\, \orcidlink{0000-0003-3478-3261}}
\address[H.~Mishra]{Department of Mathematics
University of Wisconsin-Madison
480 Lincoln Drive, Madison, WI 53706, USA.}
\email{hrishabh@cmi.ac.in}

\author[A.~Ray]{Anwesh Ray\, \orcidlink{0000-0001-6946-1559}}
\address[Ray]{Chennai Mathematical Institute, H1, SIPCOT IT Park, Kelambakkam, Siruseri, Tamil Nadu 603103, India}
\email{anwesh@cmi.ac.in}

\begin{abstract}
 Let $n$ and $k$ be positive integers, and $f_n(k)$ (resp. $g_n(k)$) be the number of unital subrings (resp. unital irreducible subrings) of $\Z^n$ of index $k$. The numbers $f_n(k)$ are coefficients of certain zeta functions of natural interest. The function $k\mapsto f_n(k)$ is multiplicative, and the study of the numbers $f_n(k)$ reduces to computing the values at prime powers $k=p^e$. Given a composition $\alpha=(\alpha_1, \dots, \alpha_{n-1})$ of $e$ into $n-1$ positive integers, let $g_\alpha(p)$ denote the number of irreducible subrings of $\Z^n$ for which the associated upper triangular matrix in Hermite normal form has diagonal $(p^{\alpha_1}, \dots, p^{\alpha_{n-1}},1)$. Via combinatorial analysis, the computation of $f_n(p^e)$ reduces to the computation of $g_\alpha(p)$ for all compositions of $i$ into $j$ parts, where $i\leq e$ and $j\leq n-1$. We extend results of Liu and Atanasov-Kaplan-Krakoff-Menzel, who explicitly compute $f_n(p^e)$ for $e\leq 8$. The case $e=9$ proves to be significantly more involved. We evaluate $f_n(p^9)$ explicitly in terms of a polynomial in n and p as a polynomial in $n$ and $p$. Our results provide further evidence for a conjecture, which states that for any fixed pair $(n,e)$, the function $p\mapsto f_n(p^e)$ is a polynomial in $p$. A conjecture of Bhargava on the asymptotics for $f_n(k)$ as a function of $k$ motivates the study of the asymptotics for $g_\alpha(p)$ for certain infinite families of compositions $\alpha$, for which we are able to obtain general estimates using techniques from the geometry of numbers.
\end{abstract}


\maketitle


\section{Introduction}\label{s 1}
\par Let $G$ be an infinite group. Given a finite index subgroup $H$, we set $[G:H]$ to denote the index of $H$ in $G$. For $k\in \Z_{\geq 1}$, let $a_k(G)$ be the number of finite-index subgroups $H$ of $G$ such that $[G:H]=k$. In \cite{grunewald1988subgroups}, Grunewald, Segal and Smith introduced the \emph{zeta function} of $G$, defined as follows
\[\zeta_G(s):=\sum_{H} [G:H]^{-s}=\sum_{k=1}^\infty a_k(G)k^{-s},\]
where the above sum runs over all finite index subgroups $H$ of $G$. The properties of such zeta functions are described in \cite{du2008zeta}. In this context, zeta functions measure \emph{subgroup growth}. For a comprehensive account of this theme, we refer to \cite{lubotzky2003subgroup}. Throughout, $n>1$ is an integer and $\Z^n$ is the ring consisting of integer $n$-tuples with componentwise addition and multiplication. It is well known that the zeta function $\zeta_{\Z^n}(s)$ is given by 
\[\zeta_{\Z^n} (s)=\zeta(s)\zeta(s-1)\dots \zeta(s-(n-1)) \] (cf. \emph{loc. cit.} for five different proofs of the above result).

\par We study distribution questions for the number of finite-index unital subrings of $\Z^n$ of prescribed index. This problem is closely related to the distributions of orders in a fixed number ring $\cO_K$, a problem studied by Brackenhoff \cite{brakenhoff2009counting}, Kaplan Marcinek and Takloo-Bighash (cf. \cite{kaplan2015distribution}). Given $k\in \Z_{\geq 1}$, following \cite{liu2007counting}, let $f_n(k)$ be the number of commutative unital subrings $S$ of $\Z^n$ such that $[\Z^n:S]=k$. The \emph{subring zeta function} is given by 
\[\zeta_{\Z^n}^R(s)=\sum_{k=1}^\infty f_n(k) k^{-s}, \]
which decomposes into an Euler product 
\[\zeta_{\Z^n}^R(s)=\prod_p \zeta_{\Z^n,p}^R(s), \] where the product is over all primes $p$. The local Euler factor at $p$ is given by 
\[\zeta_{\Z^n,p}^R(s):=\sum_{e=0}^\infty f_n(p^e) p^{-es}.\]
It is not hard to show that $\zeta_{\Z^2}^R(s)=\zeta(s)$, and that
\[\zeta_{\Z^3}^R(s)=\frac{\zeta(3s-1)\zeta(s)^3}{\zeta(2s)^2},\]
cf. \cite[Proposition 4.1]{liu2007counting}. The formula for $n=3$ was first obtained by Datskovsky and Wright \cite{datskovsky1988density}. A formula for $n=4$ was obtained by Nakagawa cf. \cite{nakagawa1996orders}. Recently, there has been interest in the following question (cf. \cite[Question 1.2]{ishamlowerbound}).
\begin{question}
Let $n, e$ be a pair of natural numbers. What can be said about $f_n(p^e)$ as a function of $p$?
\end{question}

For $e\leq 5$, Liu \cite{liu2007counting} gives explicit formulae for $f_n(p^e)$. Subsequently, these formulae were generalized by Atanasov, Kaplan, Krakoff and Menzel \cite{atanasov2021counting} for $e=6,7,8$. These computations indicate that for a fixed pair $(n,e)$, the function $p\mapsto f_n(p^e)$ is a polynomial in $p$ (cf. \cite[Question 1.13]{atanasov2021counting}), and this has been further studied by Isham in \cite{isham2022number}. The zeta function $\zeta_{\Z^n}^R(s)$ is said to be \emph{uniform} if there is a rational function $W(X,Y)\in \Q(X,Y)$ such that for every prime $p$, $\zeta_{\Z^n,p}^R(s)=W(p,p^{-s})$. It is not hard to show that the zeta function $\zeta_{\Z^n}(s)$ is uniform if for all $e$, $f_n(p^e)$ is a polynomial in $p$. For further details, we refer to \cite[section 5.1]{atanasov2021counting}. Below we summarize the known values of $f_n(p^e)$ for $e\leq 5$.
\begin{theorem}[Liu \cite{liu2007counting}]\label{theorem1.2}
With respect to notation above,
\[\begin{split}
f_n(1)= & 1,
f_n(p)=\binom{n}{2},f_n(p^2)=\binom{n}{2} + \binom{n}{3} + 3 \binom{n}{4}, \\
f_n(p^3) = & \binom{n}{2} + (p+1) \binom{n}{3} + 7 \binom{n}{4} + 10 \binom{n}{5} + 15 \binom{n}{6},\\
f_n(p^4)= & \binom{n}{2} + (3p+1) \binom{n}{3} + (p^2+p+10) \binom{n}{4}  +  (10p+21) \binom{n}{5} + 70 \binom{n}{6} \\ + & 105 \binom{n}{7} + 105 \binom{n}{8},\\
 f_n(p^5)  = & {n \choose 2} + (4p + 1){n \choose 3} + (7p^2 + p + 13){n \choose 4}
+ (p^3 + p^2 + 41p + 31){n \choose 5}\\ + & (15p^2 + 35p + 141){n \choose 6}.\\ 
\end{split}\]
\end{theorem}
For the values for $e$ in the range $6\leq e\leq 8$, we refer to \cite{kaplan2015distribution}.

\par In this manuscript, we extend the above mentioned results to $e=9$. First, we introduce some further notation.
 A subring of prime power index decomposes into a direct product of irreducible subrings, which we now proceed to describe.
\begin{definition}
A subring $L$ of $\Z^n$ of index $p^e$ is said to be \emph{irreducible} if for each vector $x=(x_1, \dots, x_n)^t\in L$, we have that 
\[x_1\equiv x_2\equiv \dots \equiv x_n\left(\op{mod} p\right).\]
\end{definition}
\noindent An irreducible subring matrix is a subring matrix in the sense of Definition \ref{subring matrix def} which is in the following Hermite normal form
\begin{equation}\label{irreducible subring matrix}
A=\begin{pmatrix}
	 p^{\alpha_1}& p a_{1,2} & p a_{1,3} & \cdots & p a_{1,n-1}& 1 \\
	& p^{\alpha_2} & p a_{2,3} & \cdots & p a_{2,n-1} & 1\\
	&& p^{\alpha_3} & \cdots & p a_{3,n-1} & 1\\
	&&& \ddots & \vdots & \vdots\\
	&&&& p^{\alpha_{n-1}} & 1\\
	&&&&& 1
\end{pmatrix},
\end{equation}
with integers $\alpha_i\ge 1$ for all $i$. Let $L$ be a subring of $\Z^n$ of index $p^{e}$ and let $v_1,\dots,v_n$ be a $\mathbb{Z}$–basis of $L$, so that each element of $L$ admits a unique expansion $v=\sum_{i=1}^n b_i v_i$ with $b_i\in\mathbb{Z}$. Writing
\[
v_i=\sum_{j=1}^n a_{i,j}\, u_j
\]
in terms of the standard basis $\{u_1,\dots,u_n\}$ of $\Z^n$, we obtain an integer matrix $A=(a_{i,j})\in M_n(\mathbb{Z})$ encoding the inclusion $L\hookrightarrow \mathbb{Z}^n$. One may choose the basis $v_1,\dots,v_n$ so that $A_L=A$ is in Hermite normal form, i.e.
\[
A_L=\begin{pmatrix}
	a_{1,1} & a_{1,2} & a_{1,3} & \cdots & a_{1,n-1}& a_{1,n} \\
	& a_{2,2} & a_{2,3} & \cdots & a_{2,n-1} & a_{2,n}\\
	&& a_{3,3} & \cdots & a_{3,n-1} & a_{3,n}\\
	&&& \ddots & \vdots& \vdots\\
	&&&& a_{n-1,n-1} & a_{n-1,n}\\
	&&&&& a_{n,n}
\end{pmatrix},
\]
with $0\le a_{i,j}<a_{i,i}$ for every $1\le i<j\le n$. Then $L$ is irreducible if and only if $A_L$ is.
\par Let $g_n(p^e)$ be the number of irreducible subrings of $\Z^n$ of index $p^e$. It is shown that every subring of $p$-power index is a direct sum of irreducible subrings, and as a consequence, we have the following recursive formula (cf. Proposition 4.4 of \cite{liu2007counting})
\begin{equation}\label{recursive formula}
    f_n(p^e)=\sum_{i=0}^e\sum_{j=1}^n {n-1\choose j-1} f_{n-j}(p^{e-i})g_j(p^i).
\end{equation}
 In order to compute $f_n(p^9)$ it thus suffices to compute $g_j(p^e)$ for all values $e\leq 9$, $1\leq j\leq n$. We remark that the convention for the definition of $g_n(p^e)$ used here is that of \cite{atanasov2021counting}, and differs from that in \cite{liu2007counting}.
 \par A composition of $m$ into $r$ parts consists of a tuple $\alpha=(\alpha_1, \dots, \alpha_r)$ of positive integers such that $\sum_{i=1}^r \alpha_i=m$. Let $\mathcal{C}_{n,e}$ be set of all compositions of $e$ into $(n-1)$ parts. Given $\alpha\in \mathcal{C}_{n,e}$, let $g_{\alpha}(p)$ be the number of irreducible subring matrices (cf. Definition \ref{subring matrix def}) of the form
 \eqref{irreducible subring matrix}.
	We find that $g_n(p^e)=\sum_{\alpha\in \mathcal{C}_{n,e}} g_\alpha(p)$, and thus to compute $g_n(p^e)$, we need to compute $g_\alpha(p)$ for all $\alpha\in \mathcal{C}_{n,e}$. Therefore, the numbers $g_\alpha(p)$ can be thought of as the basic building blocks for computing the numbers $f_n(k)$. In practice, it is possible to compute $g_\alpha(p)$ in many cases via combinatorial case by case analysis. For larger values of $e$, the computation of $f_n(p^e)$ becomes significantly more involved. The total number of new compositions that one must consider does grow exponentially, and even though one may compute the value of $g_\alpha(p)$ for various types of compositions that fit into a general framework, there are many exceptional compositions that do not fit into such a framework. Furthermore, the combinatorial (case by case) analysis does become increasingly challenging for exceptional compositions of longer length. With this in mind, we state the main result below.

 \begin{theorem}\label{main}
With respect to notation above, we have that
\begin{eqnarray*}
& & f_n(p^9) ={n \choose 2} + (p^3 + 4p^2 + 4p + 1){n \choose 3} + (11p^4 + 30p^3 + 9p^2 + p + 25){n \choose 4} \\
& & + (36p^5 + 126p^4 - 33p^3 + 92p^2 + 201p + 71){n \choose 5} \\
& & + (p^7 + 112p^6 + 133p^5 - 61p^4 + 623p^3 + 1047p^2 + 396p + 571){n \choose 6} \\
& & + (p^9 + 22p^8 + 23p^7 + 59p^6 + 31p^5 + 2032p^4\\
& & \hspace{2em} +2152p^3 + 2467p^2 + 4949p+ 2485){n \choose 7}\\
& & + (p^{10} + p^9 + 2p^8 + 2p^7 + 507p^6 + 955p^5 + 3293p^4\\
& & \hspace{2em} + 5980p^3 + 23410p^2+ 17011p + 13707){n \choose 8}\\
& & + (36p^8 + 37p^7 + 157p^6 + 1159p^5 + 8545p^4 + 34997p^3\\
& & \hspace{2em} + 59371p^2 + 93649p+ 64019){n \choose 9}\\
& & + (675p^6 + 795p^5 + 18960p^4 + 48927p^3 + 250632p^2 + 330657p+ 297103){n \choose 10}\\
& & + (990p^5 + 12540p^4 + 148830p^3 + 497640p^2 + 1157145p + 1245992){n \choose 11}\\
& & + (13860p^4 + 97020p^3 + 1049895p^2 + 2961805p + 4727041){n \choose 12}\\
& & + (135135p^3 + 1036035p^2 + 5900895p + 15346045){n \choose 13}\\
& & + (945945p^2 + 7252245p + 40500460){n \choose 14} + (4729725p + 80615535){n \choose 15}\\
& & + 110810700{n \choose 16} + 91891800{n \choose 17} + 34459425{n \choose 18}.
\end{eqnarray*}
\end{theorem}
\par We describe the method of proof in more detail. It follows from results in \cite{liu2007counting,atanasov2021counting} that $g_\alpha(p)$ can be computed for compositions $\alpha$ of one of the following types 
 \begin{itemize}
 \item $\alpha=(\beta, 1, \dots, 1)$ and $n\geq 2$, 
 \item $\alpha=(2,1\dots, 1, \beta, 1, \dots, 1)$ and $n\geq 3$. 
 \end{itemize}We refer to results in section \ref{s 2} for further details. In section \ref{s 3}, we prove some general results which allow us to compute $g_\alpha(p)$, for compositions $\alpha$ of one of the following types
 \begin{itemize}
 \item $\alpha=(\beta,1, \dots, 1, \gamma)$ is a composition of length $(n-1)\geq 3$ such that $\beta>2$ and $\gamma\geq \beta-1$,
 \item $\alpha$ is of the form $\alpha=(2,1,\dots, 1, 2, 1,\dots, 1, \beta)$, where $\beta>1$,
 \item $\alpha$ is of the form $\alpha=(2,1,\dots, 1, 2, 1,\dots, 1, \beta,1)$, where $\beta>1$,
 \item $\alpha$ is of the form $\alpha=(2,1,\dots, 1, 3, 1,\dots, 1, 2)$.
 \end{itemize}

 \par We refer to a composition that does not fit into any of the above families as an exceptional composition. In section \ref{s 4} (resp. section \ref{s 5}), we compute $g_{\alpha}(p)$ for all relevant compositions beginning with $2$ (resp. $3$). In section \ref{s 6}, we compute $g_{\alpha}(p)$ for all relevant compositions beginning with $4,5$ or $6$.
 \par Our analysis for the composition $\alpha=(3,2,2,1,1)$ leads to a finite-field point-counting problem. Let $N_p$ denote the number of solutions in $\F_p^8$ to the system
\[
\begin{split}
(x_3^2-x_3)-x_2(x_7^2-x_7)-x_1(x_5^2-x_5)&=0,\\
(x_4^2-x_4)-x_2(x_8^2-x_8)-x_1(x_6^2-x_6)&=0,\\
x_3x_4-x_2x_7x_8-x_1x_5x_6&=0.
\end{split}
\]
In Proposition \ref{prop:exceptional-point-count}, we determine the number of points on this variety exactly and prove that
\[
N_p=p^5+12p^4-20p^3+30p^2-10p.
\]
The subring matrix calculation in Theorem \ref{thm:exceptional-composition} then translates this point count into the required value of $g_\alpha(p)$. In particular, we obtain
\[
g_{(3,2,2,1,1)}(p)
=
p^7+24p^6-29p^5+21p^4-4p^3.
\]

 \par The proof of Theorem \ref{main} is provided in section \ref{s 7.1}. All known computations indicate that for any composition $\alpha$, the function $g_\alpha(p)$ is a polynomial in $p$. The conditions for a matrix $A$ give rise to polynomial conditions on the entries $a_{i,j}$, and $g_\alpha(p)$ in many cases is the number of $\F_p$-points on a scheme defined by integral polynomial equations. It is certainly of interest to note that in all the examples considered, these schemes have polynomial point count in the sense of \cite[p.616, ll. -6 to -2]{hausel2008mixed}.
 
 \par We now come to describing the general asymptotic results proved in this manuscript. The analysis of the numbers $f_n(p^e)$ can be translated into properties of the associated subring zeta functions. For instance, Isham in \cite{ishamlowerbound} proves lower bounds for $g_n(p^e)$ and deduces that the subring zeta function $\zeta_{\Z^n,p}^R(s)$ diverges for all $s$ such that $\op{Re}s\leq c_7(n)$, where $c_7(n):=\op{max}_{0\leq d \leq n-1} \frac{d(n-1-d)}{(n-1+d)}$. This comes as a consequence of proving lower bounds for $f_n(p^e)$. Some related results are also proven by Brackenhoff in \cite{brakenhoff2009counting}. Liu attributes the following conjecture regarding the asymptotics for the numbers $f_n(k)$ to Bhargava (cf. \cite[p.298]{liu2007counting}). The conjecture was communicated to Liu via personal communication, cf. \emph{loc. cit.}
 \begin{conjecture}\label{bharagava conjecture}
 For $n$ odd, $f_n(k)=O(k^{\frac{n-1}{6}+\epsilon})$, and for $n$ even, $f_n(k)=O(k^{\frac{n^2-2n}{6n-8}+\epsilon})$.
 \end{conjecture} In particular, for a fixed prime number $p$ and a fixed value of $n$, the conjecture predicts that 
 \[f_n(p^e)=\begin{cases} & O(p^{\frac{e(n-1)}{6}})\text{ if }n\text{ is odd,}\\
 & O(p^{\frac{e(n^2-2n)}{6n-8}})\text{ if }n\text{ is odd.}\\
 \end{cases}\]Conjecture \ref{bharagava conjecture} motivates our results in section \ref{s 7.2}, where we investigate the asymptotics for $g_\alpha(p)$ for certain natural families of compositions $\alpha$. Let $n > 1$ be a natural number and $t$ be in the range $1\leq t\leq n-1$. Let $\alpha = (\alpha_1,\alpha_2,\dots,\alpha_{n-1})$ be a composition of length $(n-1)$ and for $k \in \Z_{\geq 1}$, set $\alpha_{/k,t} := (\alpha_1,\alpha_2,\dots k\alpha_t\dots,\alpha_{n-2},\alpha_{n-1})$, i.e., the composition obtained upon multiplying the $t$-th coordinate of $\alpha$ by $k$. Note that $g_{\alpha_{/k,t}}(p)$ contributes to $g_n(p^{e_k})$, where $e_k:=\sum_{i=1}^{n-1}\alpha_i+(k-1)\alpha_t$. Therefore, the asymptotic growth of $g_{\alpha_{/k,t}}(p)$ as a function of $k$ provides insight into that of $g_n(p^{e_k})$. The following result is proven via a combination of techniques from the geometry of numbers and the combinatorial analysis of large subring matrices.
\begin{theorem}\label{main2}
 Let $n > 1$ be a natural number and $\alpha = (\alpha_1,\alpha_2,\dots,\alpha_{n-1})$ a composition of length $(n-1)$ such that $\alpha_i > 1$ are integers for each $i \in \{1,2,\dots,n-1\}$ and $1 \leq t\leq n-1$. Then,
 \[g_{\alpha_{/k,t}}(p) = O(p^{\gamma k(n-2)})\text{ as }k \to \infty,\]
 where $\gamma = \max\{\alpha_1,\dots,\alpha_{t}\}$. The implied constant depends on $p$ and $n$, and not on $k$.
\end{theorem}
Although the above result is weaker than the prediction of Conjecture \ref{bharagava conjecture}, the authors expect that such arguments have the potential to lead to stronger results in the future.
\subsection*{Acknowledgment} When the first draft of the article was completed, the second author was supported by the CRM-Simons postdoctoral fellowship.

\section{Preliminaries}\label{s 2}
\par Let $n\geq 1$ be an integer and $M_n(\Z)$ denote the ring of $n\times n$-matrices with integer entries. By definition, a \emph{lattice} in $\Z^n$ is a subgroup $L$ of finite index in $\Z^n$. This index is denoted $[\Z^n:L]$. Set $u_1, \dots, u_n$ to be the standard basis of $\Z^n$, with $u_i=(0,0,\dots, 0,1,0,\dots)$, with $1$ in the $i$-th position and $0$s in all other positions. Let $v_1,\dots, v_n$ be a basis of $L$, i.e., a set of vectors such that every element $v\in L$ is uniquely represented as an integral linear combination $v=\sum_{i=1}^n b_i v_i$. Expressing $v_i=\sum_{j=1}^n a_{i,j} u_j$, consider the integer matrix $A=\left(a_{i,j}\right)\in M_n(\Z)$. We may choose $v_1,\dots, v_n$ such that the associated matrix $A$ is in Hermite normal form, i.e.,
\[A=\begin{pmatrix}
	a_{1,1} & a_{1,2} & a_{1,3} & \cdots & a_{1,n-1}& a_{1,n} \\
	& a_{2,2} & a_{2,3} & \cdots & a_{2,n-1} & a_{2,n}\\
	&& a_{3,3} & \cdots & a_{3,n-1} & a_{3,n}\\
	&&& \ddots & \vdots& \vdots\\
	&&&& a_{n-1,n-1} & a_{n-1,n}\\
	&&&&& a_{n, n}
	\end{pmatrix}\]
	
	with $0\leq a_{i,j}< a_{i,i}$ for all tuples $(i,j)$ such that $1\leq i<j\leq n$. The following result is used to reinterpret counting problems for subrings of $\Z^n$ in terms of integral matrices satisfying prescribed conditions. 
	
	\begin{lemma}
	There is a bijection between lattices $L\subset \Z^n$ of index $k>0$ and integer $n\times n$ matrices $A$ in Hermite normal form with determinant $k$.
	\end{lemma}
\begin{proof}
The reader is referred to the proof of \cite[Proposition 2.1]{liu2007counting}.
\end{proof}
Thus, to a lattice $L$ we associate the matrix $A_L$ in Hermite normal form, and to a matrix $A$ in Hermite normal form, we associate a unique lattice $L_A$. 
Given two integral vectors $u=(u_1, \dots, u_n)^t$ and $v=(v_1, \dots, v_n)^t$, denote the \emph{composite} by $u\circ v:=(u_1 v_1, \dots, u_n v_n)^t$. A lattice $L$ is multiplicatively closed if $u\circ v\in L$ for all elements $u,v\in L$. A \emph{subring} of $\Z^n$ shall in this paper be taken to mean a multiplicatively closed lattice which contains the identity element $\mathbf{1}:=(1,1,\dots, 1)^t$. In particular, the index of a subring is finite. Given a positive integer $k$, let $f_n(k)$ be the number of subrings of $\Z^n$ with index equal to $k$. Define the subring zeta function as follows
\[\zeta_{\Z^n}^R(s):=\sum_{k=1}^\infty f_n(k) k^{-s}.\] The function $f_n(k)$ is multiplicative (cf. \cite[Proposition 2.7]{liu2007counting}), i.e., given two coprime integers $k_1>0$ and $k_2>0$, we have that $f_n(k_1 k_2)=f_n(k_1) f_n(k_2)$. The zeta function exhibits an Euler product 
\[\zeta_{\Z^n}^R(s)=\prod_{p}\zeta_{\Z^n,p}^R(s),\] where $\zeta_{\Z^n,p}^R(s)=\sum_{e=0}^\infty f_n(p^e)p^{-es}$. Thus, the study of the function $f_n(k)$ and the subring zeta function, comes down to the determination of $f_n(p^e)$, where $n>0, e\geq 0$ and $p$ is a prime number. More precisely, given $n$ and $e$, we would like to determine $f_n(p^e)$ as a function of $p$. 

\par We recall Liu's bijection between subrings of $\Z^n$ and integral matrices in Hermite normal form. 
\begin{proposition}\label{liu matrices}
Let $n, k>0$ be integers. The association $L\mapsto A_L$ gives a bijection between subrings $L$ of $\Z^n$ of index $k$ and matrices $A\in \op{M}_n(\Z)$ in Hermite normal form with $\det(A) = k$ and columns $v_1, \dots, v_n$ such that 
\begin{enumerate}
    \item $\mathbf{1}$ is in the column span of $A$,
    \item for all $i,j$ in the range $1\leq i, j\leq n$, $v_i\circ v_j$ is in the column span of $A$.
\end{enumerate}
\end{proposition}
\begin{proof}
The reader is referred to \cite[Proposition 2.1, 2.2]{liu2007counting}.
\end{proof}

\begin{definition}\label{subring matrix def}
A matrix $A$ satisfying the conditions of Proposition \ref{liu matrices} is referred to as a \emph{subring} matrix. The column span of $A$ is denoted $\op{Col}(A)$. 
\end{definition}

Let $\alpha=(\alpha_1, \dots, \alpha_{n-1})$ be a composition of length $(n-1)$. Recall from the introduction that $g_\alpha(p)$ is the number of irreducible subring matrices with diagonal $(p^{\alpha_1}, \dots, p^{\alpha_{n-1}}, 1)$. We recall some useful results from \cite{atanasov2021counting, liu2007counting} which shall allow us to deduce the values of $g_\alpha(p)$ for specific choices of $\alpha$. 
	
	\begin{lemma}\label{akkm lemma 3.5}
	Let $n\geq 2$ and $\alpha=(\beta, 1, \dots, 1)$ be a composition of length $(n-1)$. Then, the following assertions hold.
	\begin{enumerate}
	    \item If $\beta=2$, then $g_\alpha(p)=p^{n-2}$.
	    \item If $\beta\geq 3$, then $g_\alpha(p)=(n-1)p^{n-2}$.
	\end{enumerate}
	\end{lemma}
	\begin{proof}
	The above result is \cite[Lemma 3.5]{atanasov2021counting}.
	\end{proof}

	\begin{lemma}\label{akkm lemma 3.6}
	Let $n\geq 3$ and let $\alpha=(2,1\dots, 1, \beta, 1, \dots, 1)$ be a composition of length $(n-1)$ where $\beta$ is in the $k$-th position. Let $r=n-1-k$, i.e., the number of $1$s after $\beta$. Then, the following assertions hold. 
	\begin{enumerate}
	    \item If $\beta=2$, then, 
	    \[g_\alpha(p)=p^{n-3+r}+(r+1)p^{n-3}(p-1).\]
	    \item If $\beta\geq 3$, then \[g_\alpha(p)=(r+1)\left(p^{n-3+r}+p^{n-3}(p-1)\right).\]
	\end{enumerate}
	\end{lemma}
	\begin{proof}
	The above result is \cite[Lemma 3.6]{atanasov2021counting}.
	\end{proof}

\section{General results for computing the value of $g_\alpha(p)$}\label{s 3}
\par In this section, we prove a number of general results computing $g_\alpha(p)$ for various choices of $\alpha$. Given an irreducible subring with associated matrix $A\in M_n(\Z)$, recall that $v_1, \dots, v_n$ shall denote the columns of $A$, where $v_i$ is the $i$-th column.
\par The purpose of the lemmas in this section is to isolate families of compositions for which the subring conditions can be counted uniformly. These results account for most of the compositions that occur in the calculation of $f_n(p^9)$, and reduce the later sections to a relatively small number of exceptional compositions which require separate case-by-case arguments.
	\begin{lemma}
	 Let $n\geq 4$, $\alpha=(\beta,1, \dots, 1, \gamma)$ be a composition of length $(n-1)$ such that $\beta>2$ and $\gamma\geq \beta-1$. Then, we have 
	 \[g_\alpha(p)=p^{n-3+\lfloor \frac{\beta}{2}\rfloor}+(n-3)p^{n-2}.\]
	\end{lemma}
 \begin{proof}
 Consider a matrix $A\in M_n(\Z)$ in Hermite normal form 
  \[A=\begin{pmatrix}
	 p^{\beta}& p a_1 & pa_2 & \cdots & pa_{n-2}& 1 \\
	& p & 0 & \cdots & 0 & 1\\
	&& p & \cdots & 0 & 1\\
	&&& \ddots & \vdots& \vdots\\
	&&&& p^{\gamma} & 1\\
	&&&&& 1,
	\end{pmatrix}\] where the entries $a_1,\dots, a_{n-2}$ satisfy $0\leq a_i \leq p^{\beta-1}-1$.
	First we write down conditions for $v_j^2\in \op{Col}(A)$. Clearly, this condition is automatically satisfied for $j=1$ and $j=n$. Consider the values of $j$ that lie in the range $2\leq j\leq n-1$. We may write $j=i+1$, where $i$ lies in the range $1\leq i \leq n-2$. First, we consider the case when $i\leq n-3$. Observe that 
	\[v_{i+1}^2=(a_i^2p^2, 0, \dots, 0, p^2, 0, \dots, 0)^t.\]Note that $v_{i+1}^2$ is contained in $\op{Col}(A)$ if and only if $v_{i+1}^2-p v_{i+1}$ is contained in $\op{Col}(A)$. We find that 
	\[v_{i+1}^2-p v_{i+1}=\left((a_i^2-a_i)p^2,0, \dots,0 \right)^t.\] Therefore, $v_{i+1}^2$ is contained in $\op{Col}(A)$ if and only if 
	\begin{equation}
	   a_i^2-a_i\equiv 0\mod{p^{\beta-2}}.
	\end{equation}
	We find that $v_{n-1}^2=(a_{n-2}^2 p^2, 0, \dots, 0, p^{2\gamma}, 0)^t$ is in $\op{Col}(A)$ if and only if 
	\[v_{n-1}^2-p^\gamma v_{n-1}=\left( a_{n-2}^2p^2-a_{n-2}p^{\gamma+1}, 0,\dots, 0\right)^t\] is contained in $\op{Col}(A)$. Therefore, we deduce that $v_{n-1}^2$ is in $\op{Col}(A)$ if and only if
	\begin{equation}
	    a_{n-2}^2p^2-a_{n-2}p^{\gamma+1}\equiv 0 \mod{p^{\beta}}.
	\end{equation}
	It is assumed that $\gamma \geq \beta-1$ and hence the above condition is equivalent to 
	\[a_{n-2}^2\equiv 0\mod{p^{\beta-2}},\]
	i.e., 
	\begin{equation}
	    a_{n-2}\equiv 0 \mod{p^{\lceil \frac{\beta-2}{2} \rceil}}.
	\end{equation}
	For $1\leq i<j\leq n-2$, we find that $v_{i+1}\circ v_{j+1}$ is contained in $\op{Col}(A)$ if and only if 
	\begin{equation}
	    a_i a_j\equiv 0\mod{p^{\beta-2}}.
	\end{equation}
	Putting it all together, we find that $A$ is a subring matrix if and only if the following conditions are satisfied
	\begin{enumerate}
	    \item $a_i(a_i-1)\equiv 0\mod{p^{\beta-2}}$ for all $i$ in the range $1\leq i\leq n-3$,
	    \item $a_{n-2}\equiv 0\mod{p^{\lceil\beta/2\rceil-1}}$,
	    \item $a_i a_j\equiv 0\mod{p^{\beta-2}}$ for all $i,j$ in the range $1\leq i<j\leq n-2$.
	\end{enumerate}
	From condition (1) above, we find that at most one of $a_i\equiv 1\mod{p^{\beta-2}}$ for $i=1, \dots, n-3$. We thus are led to the following cases.
	\par \textbf{Case 1:} First, we consider the case when $a_i\equiv 0\mod{p^{\beta-2}}$ for all $i$ in the range $1\leq i\leq n-3$. Since $a_i$ is in the range $0\leq a_i<p^{\beta-1}$, we find that there are $p$ choices for each $a_i$ and $p^{\beta-\lceil \beta/2\rceil}=p^{\lfloor \beta/2\rfloor}$ choices for $a_{n-2}$. In total, we find that there are $p^{n-3+\lfloor \beta/2\rfloor}$ matrices in this case.
	\par \textbf{Case 2:} Consider the case when one of the $a_i$ satisfies $a_i\equiv 1\mod{p^{\beta-2}}$. Then, from (3) we find that \[a_{n-2}\equiv a_i a_{n-2}\equiv 0\mod{p^{\beta-2}}.\] Thus, for each index $i$ in the range $1\leq i\leq n-2$, we have $p^{n-2}$ choices for which $a_i\equiv 1\mod{p^{\beta-2}}$. There are $(n-3)$ values taken by $i$ for which $a_i\equiv 1\mod{p^{\beta-2}}$. Thus, in this case, there are $(n-3)p^{n-2}$ choices. 
	
\par Putting together the calculations from cases 1 and 2, we find that $g_\alpha(p)=p^{n-3+\lfloor \beta/2\rfloor}+(n-3)p^{n-2}$.
 \end{proof}
 
 \begin{lemma}\label{lemma 2}
 Let $\alpha$ be of the form $\alpha=(2,1,\dots, 1, 2, 1,\dots, 1, \beta)$, where $\beta>1$ and the second $2$ occurs at position $k\geq 2$, then
 \[g_\alpha(p)=p^{n-3+r}+p^{n-3}(p-1)r,\]
 where $r=n-1-k$.
 \end{lemma}
 
 \begin{proof}
 Let $A$ be an integral matrix in Hermite normal form \[A=\begin{pmatrix}
	 p^{2}& p a_1 & pa_2 & \cdots & pa_{k-1}& \cdots & \cdots & p a_{n-2} & 1 \\
	& p & 0 & \cdots & 0 &  \cdots & \cdots & 0 & 1\\
	&& p & \cdots & 0 & \cdots &  \cdots & 0  & 1\\
	&&& \ddots & \vdots& \vdots & \cdots & \vdots & \vdots\\
	&&&& p^{2} & p b_1 & \cdots & p b_r & 1\\
    &&&&  & p & \cdots & 0 & 1\\
    	&&&&& & \ddots & \vdots & 1\\
	&&&&&&& p^{\beta} & 1\\
	&&&&&&&&1.
	\end{pmatrix}\]
Since the first entry of $\alpha$ is equal to $2$, it is easy to see that $v_i^2$ is in $\op{Col}(A)$ for $i=1,\dots, k-1$. We find that $v_k^2$ is in $\op{Col}(A)$ if and only if $v_k^2-p^2 v_k$ is in $\op{Col}(A)$. Since $p^2$ divides $a_{k-1}^2p^2-p^3a_{k-1}$, this condition is seen to be satisfied. For $i=1,\dots, r$, we observe that 
\[v_{k+i}^2=(a_{k-1+i}^2 p^2, \cdots , b_i^2p^2, \cdots, p^2)^t\text{ is in }\op{Col}(A)\]if and only if 
\[v_{k+i}^2-p v_{i+k}=\left((a_{k-1+i}^2-a_{k-1+i})p^2, \dots, (b_i^2-b_i) p^2, 0,\dots,0\right)^t\] is in $\op{Col}(A)$. Subtracting $(b_i^2-b_i) v_k$ from the above, we get 
\[v_{k+i}^2-p v_{i+k}-(b_i^2-b_i) v_k=\left((a_{k-1+i}^2-a_{k-1+i})p^2-(b_i^2-b_i)a_{k-1}p, 0, \dots, 0\right)^t, \] which is in $\op{Col}(A)$ if and only if first entry is divisible by $p^2$. We have thus shown that $v_{k+i}^2$ is in $\op{Col}(A)$ if and only if 
\begin{equation}
    (b_i^2-b_i)a_{k-1}\equiv 0\mod{p}. 
\end{equation}
Now, $v_{n-1}^2=\left(a_{n-2}^2 p^2, \dots, b_r^2 p^2, 0,  \dots, 0, p^{2\beta}, 0, \dots \right)^t$ is in $\op{Col}(A)$ if and only if $v_{n-1}^2-p^\beta v_{n-1}$ is in $\op{Col}(A)$. Note that therefore, $v_{n-1}^2$ is in $\op{Col}(A)$ if and only if $b_r a_{k-1}\equiv 0\mod{p}$. 

\par It is clear that for all $1\leq i\leq k$ and all values of $j$, $v_i \circ v_j\in \op{Col}(A)$. Now consider $v_{k+i}\circ v_{k+j}$ for $1\leq i<j\leq r$, it is equal to \[\left(a_{k-1+i} a_{k-1+j}p^2, 0, \dots, 0, b_i b_j p^2, 0, \dots, 0\right)^t.\] We find that $v_{k+i}\circ v_{k+j}$ is contained in $\op{Col}(A)$ if and only if 
\[v_{k+i}\circ v_{k+j}-b_i b_j v_{k-1}\in \op{Col}(A).\] Therefore, $v_{k+i}\circ v_{k+j}$ is contained in $\op{Col}(A)$ if and only if $b_ib_ja_{k-1}\equiv 0 \mod{p}$. Therefore, $A$ is a subring matrix if and only if 
\begin{enumerate}
    \item $(b_i^2-b_i)a_{k-1}\equiv 0 \mod{p}$ for $i=1, \dots, r-1$, 
    \item $b_ra_{k-1}\equiv 0\mod{p}$, 
    \item $b_ib_ja_{k-1}\equiv 0\mod{p}$, for all $i, j$ such that $1\leq i<j\leq r$. 
\end{enumerate}
We consider two cases.
\par \textbf{Case 1:} Consider the case when $a_{k-1}=0$. In this case, the total number of choices is the total number of choices of $(a_1, \dots, a_{k-2}, a_{k}, \dots, a_{n-2})$ and $(b_1, \dots, b_r)$. Thus, the total number of choices are $p^{n-3+r}$.
\par \textbf{Case 2:} Consider the other case, i.e., when $a_{k-1}\neq 0$. Since $a_{k-1}<p$, we find that $p\nmid a_{k-1}$. Therefore, the conditions are as follows
\begin{enumerate}
    \item $b_i^2-b_i\equiv 0\mod{p}$, 
    \item $b_r\equiv 0\mod{p}$, 
    \item $b_i b_j\equiv 0\mod{p}$, $1\leq i<j\leq r$.
\end{enumerate}
Note that all elements $b_i$ satisfy the bounds $0\leq b_i<p$. Therefore, at most one of the $b_i$ satisfies $b_i=1$. Subdivide into two cases, first consider the case when all the $b_i$ are equal to $0$ and then consider the case when at most one of the $b_i$ is equal to $1$. The total number of choices is $p^{n-3}(p-1)r$. Upon adding up all the conclusions from the case decomposition above, we find that
\[g_\alpha(p)=p^{n-3+r}+p^{n-3}(p-1)r.\]
  \end{proof}
 
\begin{lemma}\label{lemma 3}
 Let $\alpha$ be of the form $\alpha=(2,1,\dots, 1, 2, 1,\dots, 1, \beta,1)$, where $\beta>1$, the second $2$ occurs at position $k\geq 2$ and set $r:=n-1-k$. Then, the following assertions hold.
 \begin{enumerate}
    \item If $\beta = 2$, then we find that
    \[g_\alpha(p)=2p^{n-3+r} + p^{n-5+r}(p-2) + 2rp^{n-3}(p-1) + p^{n-4}(p-1)(p-2)(p+r-2).\]
    \item If $\beta \geq 3$, then,
    \[g_\alpha(p)=2p^{n-3+r} + 2p^{n-5+r}(p-1) + 2rp^{n-3}(p-1) + 2p^{n-4}(p-1)^2(p+r-2).\]
\end{enumerate}
\end{lemma}

\begin{proof}
 Let $A$ be an integral matrix in Hermite normal form \[A=\begin{pmatrix}
	 p^{2}& p a_1 & pa_2 & \cdots & pa_{k-1}& \cdots & \cdots &\cdots& p a_{n-2} & 1 \\
	& p & 0 & \cdots & 0 &  \cdots & \cdots &\cdots& 0 & 1\\
	&& p & \cdots & 0 & \cdots &  \cdots &\cdots& 0  & 1\\
	&&& \ddots & \vdots& \vdots & \cdots &\cdots& \vdots & \vdots\\
	&&&& p^{2} & p b_1 & \cdots &\cdots& p b_r & 1\\
    &&&&  & p & \cdots &\cdots& 0 & 1\\
    	&&&&& & \ddots &\cdots& \vdots & 1\\
	&&&&&&& p^{\beta} & p c_1 & 1\\
	&&&&&&&&p&1\\
	&&&&&&&&&1.
	\end{pmatrix}\]
Since the first entry is $p^2$, it is clear that $v_i^2$ is contained in $\op{Col}(A)$ for $i=1,2,\dots,k$. For $i = 1,2,\dots,r-2$, we find that $v_{k+i}^2$ is contained in $\op{Col}(A)$ if and only if $v_{k+i}^2-pv_{k+i}$ is contained in $\op{Col}(A)$. So,
\[(a_{k-1+i}^2p^2-a_{k-1}p^2,\dots,b_i^2p^2-b_ip^2,\dots)^t\]
should be in $\op{Col}(A)$. It is clear that this is true if and only if
\begin{equation}
    a_{k-1}(b_i^2-b_i) \equiv 0 \mod{p}
\end{equation}
Similarly we see that $v_{n-2}^2 \in \op{Col}(A)$ if and only if $v_{n-2}^2 - p^{\beta}v_{n-2} \in \op{Col}(A)$. Now,
\[(a_{n-3}^2p^2-a_{n-3}p^{\beta+1},\dots,b_{r-1}^2p^2-b_{r-1}p^{\beta+1},\dots)^t\]
should be in $\op{Col}(A)$. It is again clear that this is true if and only if
\begin{equation}
    a_{k-1}b_{r-1} \equiv 0 \mod p
\end{equation}
as $\beta > 1$. Again, $v_{n-1}^2 \in \op{Col}(A)$ if and only if $v_{n-1}^2 - pv_{n-1} \in \op{Col}(A)$, so
\[(a_{n-2}^2p^2-a_{n-2}p^2,\dots,(b_r^2-b_r)p^2,\dots,(c_1^2-c_1)p^2, 0, 0)^t\]
should be in $\op{Col}(A)$. Then we have $c_1^2-c_1 \equiv 0 \mod p^{\beta-2}$, and
\begin{equation}
 b_{r-1}(c_1^2-c_1) \equiv 0 \mod p^{\beta-1}
 \end{equation}
 \begin{equation}
 a_{k-1}\left(b_r^2-b_r-\frac{b_{r-1}(c_1^2-c_1)}{p^{\beta-1}}\right)+a_{n-3}\frac{c_1^2-c_1}{p^{\beta-2}} \equiv 0 \mod p
\end{equation}
The requirement that $v_i\circ v_j$ belongs to $\op{Col}(A)$ for $i\neq j$ translates the following conditions on entries of the matrix (for all $1 \leq i < j \leq r$, and $\{i,j\} \neq \{r-1,r\}$)
\begin{equation*}
    a_{k-1}b_ib_j \equiv 0 \mod p
\end{equation*}
First, we suppose that $\beta \geq 3$, then we have conditions,
\begin{enumerate}
    \item $a_{k-1}(b_i^2-b_i) \equiv 0 \mod{p}$, $i=1,\dots,r-2$,
    \item $b_{r-1}a_{k-1} \equiv 0 \mod p$,
    \item $c_1^2-c_1 \equiv 0 \mod p^{\beta-2}$,
    \item $b_{r-1}(c_1^2-c_1) \equiv 0 \mod p^{\beta-1}$,
    \item $a_{k-1}\left(b_r^2-b_r-\frac{b_{r-1}(c_1^2-c_1)}{p^{\beta-1}}\right)+a_{n-3}\frac{c_1^2-c_1}{p^{\beta-2}} \equiv 0 \mod p$,
    \item $a_{k-1}b_ib_j \equiv 0 \mod p$, $1 \leq i < j \leq r$, and $\{i,j\} \neq \{r-1,r\}$.
\end{enumerate}
We consider two cases.
\par \textbf{Case 1 :} First, we consider the case when $a_{k-1} = 0$. Then, the above equations reduce to the following:
\begin{enumerate}
    \item $c_1^2-c_1 \equiv 0 \mod p^{\beta-2}$,
    \item $b_{r-1}(c_1^2-c_1) \equiv 0 \mod p^{\beta-1}$,
    \item $\frac{a_{n-3}(c_1^2-c_1)}{p^{\beta-2}} \equiv 0 \mod p$.
\end{enumerate}
If $c_1 \in \{0,1\}$, then the number of such matrices is $2p^{n-3+r}$. On the other hand, if $c_1 \notin \{0,1\}$, then, there are $2p-2$ choices for $c_1$. Since $b_{r-1} = 0, a_{n-3}=0$, we deduce that there are $2(p-1)p^{n+r-5}$ more matrices.
\par \textbf{Case 2 :} Suppose $a_{k-1} \neq 0$, then equations reduce to
\begin{enumerate}
\item $b_i^2-b_i \equiv 0 \mod p$, $i =1, \dots, r-2$,
\item $b_{r-1}=0$,
\item $c_1^2-c_1 \equiv 0 \mod p^{\beta-2}$,
\item $a_{k-1}\left(b_r^2-b_r\right)+a_{n-3}\frac{c_1^2-c_1}{p^{\beta-2}} \equiv 0 \mod p$,
\item $b_ib_j \equiv 0 \mod p$, $1 \leq i < j \leq r$, and $\{i,j\} \neq \{r-1,r\}$
\end{enumerate}
We divide into two sub-cases.
\begin{itemize}
    \item \textbf{Case 2A :} Suppose that $c_1 \in \{0,1\}$, then we get that $b_r \in \{0,1\}$ and as in Lemma \ref{lemma 2}, we get that at most one of the $b_i$ satisfies $b_i=1$. Subdivide into two cases, first consider the case when all the $b_i$ are equal to $0$ and then consider the case when at most one of the $b_i$ is equal to $1$. Therefore, the total number of choices is $2rp^{n-3}(p-1)$.
    \item \textbf{Case 2B :} Consider the other subcase, i.e., $c_1 \notin \{0,1\}$. Then, there is a unique value of $a_{n-3}$ for each value of $b_r, a_{k-1}, c_1$, so we get $2p^{n-3}(p-1)^2(p+r-2)$ matrices
\end{itemize}
Putting it all together, the result is proven. The case $\beta = 2$ is similar, and the number of matrices change only in Case 1, with $c_1\notin \{0,1\}$ and Case 2B.
\end{proof}

 \begin{lemma}\label{lemma 4}
 Let $\alpha$ be of the form $\alpha=(2,1,\dots, 1, 3, 1,\dots, 1, 2)$, where $3$ occurs at position $k\geq 2$. Then, we have that
 \[g_{\alpha}(p) = rp^{n-3+r} + p^{n-3}(p-1)(p+r-1),\]
 where $r:=n-1-k$.
 \end{lemma}
 \begin{proof}
   Let $A$ be an integral matrix in Hermite normal form \[A=\begin{pmatrix}
	 p^{2}& p a_1 & pa_2 & \cdots & pa_{k-1}& \cdots & \cdots & p a_{n-2} & 1 \\
	& p & 0 & \cdots & 0 &  \cdots & \cdots & 0 & 1\\
	&& p & \cdots & 0 & \cdots &  \cdots & 0  & 1\\
	&&& \ddots & \vdots& \vdots & \cdots & \vdots & \vdots\\
	&&&& p^{3} & p b_1 & \cdots & p b_r & 1\\
    &&&&  & p & \cdots & 0 & 1\\
    	&&&&& & \ddots & \vdots & 1\\
	&&&&&&& p^{2} & 1\\
	&&&&&&&&1.
	\end{pmatrix}\]
	Since the first diagonal entry is $p^2$, we see that $v_i^2$ is contained in $\op{Col}(A)$ for all $i=1,\dots,k$. Note that $v_{k+i}^2 \in \op{Col}(A)$ if and only if $v_{k+i}^2-pv_{k+i} \in \op{Col}(A)$ for $i = 1,2,\dots,r-1$. Hence, we find that
	\[\left((a_{k-1+i}^2-a_{k-1+i})p^2,0,\dots,0,(b_i^2-b_i)p^2,0,\dots,0\right)^t\]
	should be in $\op{Col}(A)$. We deduce that that this is the case if and only if the following conditions are satisfied
	\begin{equation}
	\begin{split}
	   & b_i^2-b_i \equiv 0 \mod p,\\
	   &  a_{k-1}(b_i^2-b_i) \equiv 0 \mod p^2.
	 \end{split}
	\end{equation}
	It is also required that $v_{n-1}^2$ is contained in $\op{Col}(A)$, which is the case if and only if $v_{n-1}^2 - p^2v_{n-1}$ is contained in $ \op{Col}(A)$. Therefore, we find that
	\[\left((a_{n-2}^2-a_{n-2}p)p^2,0,\dots,0,(b_r^2-b_r p)p^2,0,\dots,0\right)^t\]
	should be in $\op{Col}(A)$. We get that $b_r^2 - b_r p \equiv 0 \mod p$ and $a_{k-1}(b_r^2-b_r p) \equiv 0 \mod p^2$, therefore, we find that $p \mid b_r$. Using similar arguments we deduce that $v_{k+i}v_{k+j} \in \op{Col}(A)$ for $1 \leq i < j \leq r$ if and only if the following conditions are satisfied
	\begin{equation}
 \begin{split}
	    & b_i b_j \equiv 0 \mod p,\\
	    & a_{k-1}b_i b_j \equiv 0 \mod p^2.
     \end{split}
	\end{equation}
It is clear that $v_i\circ v_j$ is contained in $\op{Col}(A)$ for $1 \leq i < j \leq k$. Therefore, we have following conditions on the entries of $A$
	\begin{enumerate}
	    \item $b_i^2-b_i \equiv 0 \mod p$ and $a_{k-1}(b_i^2-b_i) \equiv 0 \mod p^2$ for $i=1,\dots,r-1$,
	    \item $b_r = b_r' p$, where $0 \leq b_r' \leq p-1$,
	    \item $b_i b_j \equiv 0 \mod p$ and $a_{k-1}b_i b_j \equiv 0 \mod p^2$ for all values $1 \leq i < j \leq r$.
	\end{enumerate}
	We consider two cases as follows.
	\par \textbf{Case 1 :} First, we consider the case when $a_{k-1} = 0$. In this case, the conditions on $A$ reduce to the following
	\begin{enumerate}
	    \item $b_i^2-b_i \equiv 0 \mod p$ for all $i=1,\dots,r-1$,
	    \item $b_r = b_r' p$, where $0 \leq b_r' \leq p-1$,
	    \item $b_i b_j \equiv 0 \mod p$ for $1 \leq i < j \leq r$.
	\end{enumerate}
	As in Lemma \ref{lemma 3}, either all the values of $b_i$ for $1 \leq i \leq r-1$ are $0 \mod p$, or at most one of these values is $1 \mod p$. Further dividing into cases, it is easy to see that the number of matrices is $rp^{n-3+r}$.
	\par \textbf{Case 2 :} Next, we consider the other case, namely assume that $a_{k-1} \neq 0$. The conditions on $A$ then reduce to the following
	\begin{enumerate}
	    \item $b_i^2-b_i \equiv 0 \mod p^2$ for all $i=1,\dots,r-1$,
	    \item $b_r = b_r' p$, where $0 \leq b_r' \leq p-1$,
	    \item $b_i b_j \equiv 0 \mod p^2$ for all $(i,j)$ satisfying $1 \leq i < j \leq r$.
	\end{enumerate}
	Consequently, it follows that either $b_i=0$ for all $i=1,\dots,r-1$, or at most one of them is $1$. Therefore, we find that the number of matrices in this case is equal to $p^{n-3}(p-1)(p+r-1)$ matrices. Adding up our conclusions, we prove the assertion of the lemma. 
 \end{proof}

\section{Calculating the values of $g_\alpha(p)$ for compositions beginning with $2$}\label{s 4}
\subsection{Compositions of length $4$} We consider the compositions $\alpha$ of length $4$. In all, there are $15$ of them, listed below
\begin{table}[H]
\centering
 \begin{tabular}{c c } 
 \hline
 $(2,5,1,1)$ & $(2,1,2,4)$\\ 
 $(2,1,5,1)$ & $(2,3,3,1)$\\
 $(2,1,1,5)$ & $(2,3,1,3)$\\
 $(2,4,2,1)$ & $(2,1,3,3)$ \\
 $(2,4,1,2)$ & $(2,3,2,2)$ \\
 $(2,2,4,1)$ & $(2,2,3,2)$\\
 $(2,2,1,4)$ & $(2,2,2,3)$.\\
 $(2,1,4,2)$ & \\
\hline
\end{tabular}
\end{table}
\noindent It follows directly from Lemmas \ref{akkm lemma 3.6}, \ref{lemma 2} and \ref{lemma 3} that
\begin{table}[H]
\centering
 \begin{tabular}{c c } 
 \hline
 $g_{(2,5,1,1)}(p)$ & $3p^4+3p^3-3p^2$\\ 
 $g_{(2,1,5,1)}(p)$ & $4p^3-2p^2$\\
 $g_{(2,1,1,5)}(p)$ & $p^3$\\
 $g_{(2,2,4,1)}(p)$ & $4p^4+2p^3-4p^2$\\
 $g_{(2,2,1,4)}(p)$ & $p^4+2p^3-2p^2$\\
 $g_{(2,1,2,4)}(p)$ & $2p^3-p^2$.\\
\hline
\end{tabular}
\end{table}
\noindent We compute $g_\alpha(p)$ for the rest of the compositions. We note that many arguments are similar, and we summarize the arguments that tend to repeat.  
\begin{lemma}\label{lemma 4.1}
With respect to notation above, we have that $g_{(2,1,4,2)}(p) =3p^4-2p^3$.
\end{lemma}
\begin{proof} Let $A$ be the matrix,
$$\begin{pmatrix}p^2 & a_1p & a_2p & a_3p & 1\\ 0 & p & 0 & 0 & 1\\
0 & 0 & p^4 & b_1p & 1\\
0 & 0 & 0 & p^2 & 1\\
0 & 0 & 0 & 0 & 1\end{pmatrix}$$
where $0 \leq a_1 \leq p-1$ and $0 \leq b_1 \leq p^3-1$.
We arrive at conditions for $A$ to be a subring matrix. It is easy to see that
\begin{itemize}
    \item $v_2^2$,$v_3^2$,$v_2\circ v_3$, $v_2\circ v_4$, $v_3\circ v_4$ are in $\op{Col}(A)$.
    \item We find that $v_4^2$ is in $\op{Col}(A)$ if and only if $v_4^2-p^2v_4$ is contained in $\op{Col}(A)$. In other words,
    \[\begin{split}& (a_3^2p^2,
    0,
    b_1^2p^2,
    p^4,
    0)^t - p^2(a_3p,
    0,
    b_1p,
    p^2,
    0)^t \\ =& \left((a_3^2p^2-a_3p^3), 0, (b_1^2p^2-b_1p^3),0,0\right)^t\end{split}\] is in $\op{Col}(A)$. The expression $(b_1^2p^2-b_1p^3)=b_1(b_1-p)p^2$ must be divisible by $p^4$ and moreover, we find that $v_4^2-p^2v_4$ is contained in $\op{Col}(A)$ if and only if $v_4^2-p^2v_4-\frac{b_1(b_1-p)}{p^2}v_3$ is in $\op{Col}(A)$. This holds if and only if in addition, we have that $b_1(b_1-p)a_2 \equiv 0 \mod p^3$.

\end{itemize}
We deduce from above that the necessary conditions are as follows
\begin{enumerate}
    \item\label{lemma 4.1 e1} $b_1(b_1-p) \equiv 0 \mod p^2$,
    \item\label{lemma 4.1 e2} $b_1(b_1-p)a_2 \equiv 0 \mod p^3$.
\end{enumerate}
From equation \eqref{lemma 4.1 e1}, we get that $b_1 = b_1'p$, where $0 \leq b_1' \leq p^2-1$. Equation \eqref{lemma 4.1 e2} asserts that $b_1'(b_1'-1)a_2 \equiv 0 \mod p$. Consider the following case decomposition.
\begin{itemize}
    \item \textbf{Case 1 :} Assume that $a_2 = 0$. In this case, the number of such matrices is $p^4$.
    \item \textbf{Case 2 :} Consider the other case, i.e., $a_2\neq 0$. Then, we have $b_1'(b_1'-1) \equiv 0 \mod p$. Hence, the total number of such matrices is easily seen to be $2p^3(p-1)$.
\end{itemize}
We conclude from the above that $g_{(2,1,4,2)}(p) = p^4+2p^3(p-1)=3p^4-2p^3$.
\end{proof}
\begin{lemma}\label{lemma 4.2}
With respect to notation above, we find that $g_{(2,4,1,2)}(p) = p^5 + 3p^4 - p^3 - p^2$.
\end{lemma}
\begin{proof} 
Let $A$ be the matrix
$$\begin{pmatrix}p^2 & a_1p & a_2p & a_3p & 1\\ 0 & p^4 & b_1p & b_2p & 1\\
0 & 0 & p & 0 & 1\\
0 & 0 & 0 & p^2 & 1\\
0 & 0 & 0 & 0 & 1\end{pmatrix}$$
where $0 \leq a_i \leq p-1$ and $0 \leq b_j \leq p^3-1$.
We obtain conditions for $A$ to be a subring matrix. By the same arguments as in Lemma \ref{lemma 4.1}, we find that
\begin{itemize}
    \item $v_2^2$, $v_2\circ v_3$, $v_2\circ v_4$ are in $\op{Col}(A)$.
    \item We find that $v_3^2$ is in $\op{Col}(A)$ if and only if $b_1(b_1-1) \equiv 0 \mod p^2$ and $b_1(b_1-1)a_1 \equiv 0 \mod p^3$.
    \item We find that $v_4^2$ is in $\op{Col}(A)$ if and only if
    $b_2 = b_2'p$, $0 \leq b_2' \leq p^2-1$ and $b_2'(b_2'-1)a_1 \equiv 0 \mod p$.
    \item We find that $v_3\circ v_4$ is in $\op{Col}(A)$ if and only if $b_1b_2' \equiv 0 \mod p$ and $a_1b_1b_2' \equiv 0 \mod p^2$.
\end{itemize}
Hence, $A$ is a subring matrix if and only if
\begin{enumerate}
    \item $b_1(b_1-1) \equiv 0 \mod p^2$,
    \item $a_1b_1(b_1-1) \equiv 0 \mod p^3$,
    \item $b_2 = b_2'p$ for $0 \leq b_2' \leq p^2-1$,
    \item $a_1b_2'(b_2'-1) \equiv 0 \mod p$,
    \item $b_1b_2' \equiv 0 \mod p$,
    \item $a_1b_1b_2' \equiv 0 \mod p^2$.
\end{enumerate}
In order to compute the total number of matrices satisfying all of the above conditions, we consider the following cases.
\begin{itemize}
    \item \textbf{Case 1 :} Assume that $a_1 = 0$. If $b_1 \equiv 0 \mod p^2$ then number of such matrices is $p^5$ otherwise $b_1 \equiv 1 \mod p^2$, number of such matrices is $p^4$.
    \item \textbf{Case 2 :} Consider the case when $a_1 \neq 0$. We get $b_1 = 0$ or $b_1 = 1$. If $b_1 = 0$, number of such matrices is $2p^3(p-1)$. Otherwise $b_1 = 1$ and number of such matrices is $p^2(p-1)$ (as $b_2 = 0$ in this case).
\end{itemize}
Therefore, we find that $g_{(2,4,1,2)}(p) = p^5+p^4+2p^3(p-1)+p^2(p-1)= p^5 + 3p^4 - p^3 - p^2$.
\end{proof}
\begin{lemma}
We have that \[\begin{split}& g_{(2,3,1,3)}(p) = 3p^4-p^2, \\
& g_{(2,1,3,3)}(p) = p^4,\\
& g_{(2,2,3,2)}(p)=p^5+p^4-p^3, \\
& g_{(2,2,2,3)}(p)=2p^4-p^2.\end{split}\]
\end{lemma}
\begin{proof}
The proof is similar to Lemma \ref{lemma 4.2}, and we omit it.
\end{proof}
\begin{lemma}\label{lemma 4.4}
We have that $g_{(2,3,2,2)}(p)=p^5+3p^4-5p^3+2p^2$.
\end{lemma}
\begin{proof}
Let $A$ be the matrix
$$\begin{pmatrix}p^2 & a_1p & a_2p & a_3p & 1\\ 0 & p^3 & b_1p & b_2p & 1\\
0 & 0 & p^2 & c_1p & 1\\
0 & 0 & 0 & p^2 & 1\\
0 & 0 & 0 & 0 & 1\end{pmatrix}$$
where $0\leq a_i\leq p-1$, $0\leq b_j\leq p^2-1$ and $0\leq c_1\leq p-1$. As in the preceding lemmas, we first impose the conditions that the squares and pairwise composites of the columns lie in $\op{Col}(A)$. From $v_3^2\in\op{Col}(A)$ we obtain
\[b_1=b_1'p\quad\text{where}\quad 0\leq b_1'\leq p-1.\]
The condition coming from $v_4^2$ is
\begin{equation}\label{lemma 4.4 corrected 1}
b_2^2-c_1^2b_1'\equiv 0\mod p,
\end{equation}
and, after dividing the corresponding first-coordinate relation by $p$, we obtain
\begin{equation}\label{lemma 4.4 corrected 2}
\left(\frac{b_2^2-c_1^2b_1'}p-b_2\right)a_1+c_1(c_1-1)a_2\equiv0\mod p.
\end{equation}
Finally, the condition $v_3\circ v_4\in\op{Col}(A)$ is
\begin{equation}\label{lemma 4.4 corrected 3}
b_1'a_1(b_2-c_1)\equiv0\mod p.
\end{equation}
These are the only additional conditions. Counting the solutions of \eqref{lemma 4.4 corrected 1}--\eqref{lemma 4.4 corrected 3}, first according as $a_1=0$ or $a_1\neq0$ and then according as $b_1'=0$ or $b_1'\neq0$, gives respectively
\[p^4+2p^3-p^2\quad\text{and}\quad 2p^4-7p^3+3p^2.\]
Adding the two contributions gives
\[g_{(2,3,2,2)}(p)=p^5+3p^4-5p^3+2p^2.\]
\end{proof}
\begin{lemma} We have that \[\begin{split} & g_{(2,4,2,1)}(p) = 2p^5+12p^4-17p^3+4p^2, \\ 
& g_{(2,3,3,1)}(p) = 12p^4-10p^3+2p^2.\end{split}\]
\end{lemma}
\begin{proof}
The proof is similar to that of Lemma \ref{lemma 4.4}, and is omitted.
\end{proof}

\subsection{Compositions of length $5$} 
In this subsection we consider the compositions of length $5$ that begin with $2$. We need to compute  $g_{\alpha}(p)$ for $20$ compositions of this form. They are listed below,\\
\begin{table}[hbt!]
\centering
\begin{tabular}{c c } 
\hline
$(2,4,1,1,1)$ & $(2,1,4,1,1)$\\ 
$(2,1,1,4,1)$ & $(2,1,1,1,4)$\\
$(2,3,2,1,1)$ & $(2,3,1,2,1)$\\
$(2,3,1,1,2)$ & $(2,2,3,1,1)$\\
$(2,2,1,3,1)$ & $(2,2,1,1,3)$\\
$(2,1,3,1,2)$ & $(2,1,3,2,1)$\\
$(2,1,2,3,1)$ & $(2,1,2,1,3)$\\
$(2,1,1,3,2)$ & $(2,1,1,2,3)$\\
$(2,2,2,2,1)$ & $(2,2,2,1,2)$\\
$(2,2,1,2,2)$ & $(2,1,2,2,2)$\\
\hline
\end{tabular}
\caption{Compositions of length $5$ that begin with $2$}
\end{table}\\
In the next two lemmas we obtain values of $g_{\alpha}(p)$ for all the compositions listed above.
\begin{lemma}
We have the following values
\begin{table}[H]
\centering
 \begin{tabular}{c c } 
 \hline
 
 $g_{(2,4,1,1,1)}(p)$ & $4p^6 + 4p^4 - 4p^3$\\ 
 $g_{(2,1,4,1,1)}(p)$ & $3p^5+3p^4-3p^3$\\
 $g_{(2,1,1,4,1)}(p)$ & $4p^4-2p^3$\\
 $g_{(2,1,1,1,4)}(p)$ & $p^4$ \\
 $g_{(2,3,1,1,2)}(p)$ & $3p^6+p^5+p^4-2p^3$ \\
 $g_{(2,1,3,1,2)}(p)$ & $3p^5-p^3$\\
 $g_{(2,1,1,3,2)}(p)$ & $p^5$\\
 $g_{(2,2,1,1,3)}(p)$ & $p^6+3p^4-3p^3$\\
 $g_{(2,1,2,1,3)}(p)$ & $p^5+2p^4-2p^3$\\
 $g_{(2,1,1,2,3)}(p)$ & $2p^4-p^3$\\
 $g_{(2,2,1,3,1)}(p)$ & $2p^6+4p^5+2p^4-8p^3+2p^2$\\
 $g_{(2,1,2,3,1)}(p)$ & $4p^5+2p^4-4p^3$\\
 
\hline
\end{tabular}
\end{table}
\end{lemma}
\begin{proof}
    The result follows immediately from Lemmas \ref{akkm lemma 3.6}, \ref{lemma 2}, \ref{lemma 3} and \ref{lemma 4}.
\end{proof}
Next, we obtain values of $g_{\alpha}(p)$ for remaining $8$ compositions.
\begin{lemma}
We have the following values,
\begin{table}[H]
\centering
 \begin{tabular}{c c} 
 \hline
 $g_{(2,3,2,1,1)}(p)$ & $16p^6+p^5-31p^4+27p^3-6p^2$\\ 
 $g_{(2,2,3,1,1)}(p)$ & $6p^6+15p^5-27p^4+12p^3-3p^2$\\ 
 $g_{(2,1,3,2,1)}(p)$ & $10 p^5 - 11 p^4 + 4 p^3$\\
 $g_{(2,3,1,2,1)}(p)$ & $7p^6+4p^5-8p^4-p^3+2p^2$\\
 $g_{(2,2,2,2,1)}(p)$ & $8p^6+5p^5-20p^4+12p^3-4p^2$\\
 $g_{(2,2,2,1,2)}(p)$ & $3p^6+5p^5-11p^4+6p^3-2p^2$\\
 $g_{(2,2,1,2,2)}(p)$ & $p^6+2p^5-3p^3+p^2$\\
 $g_{(2,1,2,2,2)}(p)$ & $2 p^5 - p^3$\\

\hline
\end{tabular}
\end{table}
\end{lemma}
\begin{proof}We will prove this for the composition $(2,1,2,2,2)$, expression for other compositions can be obtained in a similar way. Let $A$ be a matrix in Hermite normal form of following type,
$$
    \begin{pmatrix}
        p^2 & a_1p & a_2p & a_3p & a_4p & 1\\
         & p & 0 & 0 &  0 & 1\\
          & & p^2 & b_1p & b_2p & 1\\
          &&& p^2 & c_1p & 1\\
          &&&& p^2& 1\\
          &&&&&1\\
    \end{pmatrix}
    $$
    First, we determine the conditions on entries of the matrix which make $A$ a subring matrix.
    \begin{enumerate}
        \item First, we note that $v_2^2 \in \Col(A)$ and $v_3^2 \in \Col(A)$.
        \item Next, we know that $v_4^2 \in \Col(A)$ if and only $v_4^2-p^2v_4 \in \Col(A)$. Now, this is the case if and only if $b_1a_2 \equiv 0 \mod p$.
        \item Similarly, $v_5^2 \in \Col(A)$ if and only if $v_5^2-p^2v_5 \in \Col(A)$ and this is true if and only if $b_1c_1 \equiv 0 \mod p$ and $b_2^2a_2 - c_1^2a_3 \equiv 0 \mod p$.
        \item We also note that if above conditions are satisfied then $v_2\circ v_3, v_2\circ v_4, v_2\circ v_5, v_3\circ v_4, v_3\circ v_5$ and $v_4\circ v_5$ are in $\Col(A)$.
    \end{enumerate}
    The conditions we get are,
    \begin{itemize}
        \item $b_1a_2 = 0$,
        \item $b_1c_1 = 0$,
        \item $b_2^2a_2-c_1^2a_3 \equiv 0 \mod p$.
    \end{itemize}
    First we assume that $b_1 = 0$. Now further if $c_1 = 0$ then we get $p^3(2p-1)$ matrices. If $c_1 \neq 0$ then we get $p^4(p-1)$ matrices. Next, we assume that $b_1 \neq 0$ then $a_2 = 0$ and $c_1=0$, therefore we get $p^4(p-1)$ matrices. Therefore, we get that $g_{(2,1,2,2,2)}(p) = 2 p^5 - p^3$.
\end{proof}
  
\section{Calculating the values of $g_\alpha(p)$ for compositions beginning with $3$}\label{s 5}
\subsection{Compositions of length $4$} 
\noindent In this section, we compute $g_{\alpha}(p)$ for $10$ compositions $\alpha$ beginning with $3$ and of length $4$. They are listed as follows
\begin{table}[H]
\centering
 \begin{tabular}{c c } 
 \hline
 $(3,4,1,1)$ & $(3,3,1,2)$\\ 
 $(3,1,4,1)$ & $(3,2,3,1)$\\
 $(3,1,1,4)$ & $(3,1,2,3)$\\
 $(3,2,1,3)$ & $(3,3,2,1)$ \\
 $(3,1,3,2)$ & $(3,2,2,2)$.\\
\hline
\end{tabular}
\end{table}
We make note of some known computations.
\begin{lemma}
The following values of $g_\alpha(p)$ are known
\begin{table}[H]
\centering
 \begin{tabular}{c c } 
 \hline
 $g_{(3,4,1,1)}(p)$ & $15p^4-9p^3+3p^2$\\ 
 $g_{(3,1,4,1)}(p)$ & $2p^4+6p^3-2p^2$\\
 $g_{(3,1,1,4)}(p)$ & $3p^3$.\\
\hline
\end{tabular}
\end{table}
\end{lemma}
\begin{proof}
The above result follows from \cite[Lemma 4.6]{isham2022number}.
\end{proof}
\begin{lemma}\label{(3,2,1,3)}
We have that $g_{(3,2,1,3)}(p) = 5p^4-4p^3+2p^2.$
\end{lemma}
\begin{proof}
    Let $A$ be any integer matrix in Hermite normal form of the type
    $$
    \begin{pmatrix}
        p^3 & a_1p & a_2p & a_3p & 1\\
         & p^2 & b_1p & b_2p & 1\\
          & & p & 0 & 1\\
          &&& p^3 & 1\\
          &&&& 1\\
    \end{pmatrix}
    $$
    We determine the conditions so that $A$ is a subring matrix.
    \begin{enumerate}
        \item It is required that $v_2^2 \in \Col(A)$ and this is the case if and only if
        \[(a_1^2p^2, p^4, 0, 0, 0)^t - p^2(a_1p, p^2, 0, 0, 0)^t\]
        is in $\Col(A)$. We deduce that $a_1 \equiv 0 \mod p$, hence $a_1 = a_1'p$ where $0 \leq a_1' \leq p-1$.
        \item It is required that $v_3^2 \in \Col(A)$; clearly this is equivalent to the statement that $v_3^2 - pv_3 \in \Col(A)$. Now, it is easy to see that latter condition is true if and only if
        \begin{equation}\label{(3,2,1,3)1}
            (a_2^2 - a_2) - a_1'(b_1^2-b_1) \equiv 0 \mod p.
        \end{equation}
        \item Similarly, it is clear that $v_4^2 \in \Col(A)$ if and only if $v_4^2 - p^3v_4 \in \Col(A)$ and latter condition is equivalent to
        \begin{equation}\label{(3,2,1,3)2}
            a_3^2-b_2^2a_1' \equiv 0 \mod p.
        \end{equation}
        \item One can easily check that $v_2\circ v_3$ and $v_2\circ v_4$ are both in $\Col(A)$ and $v_3\circ v_4 \in \Col(A)$ if and only
        \begin{equation}\label{(3,2,1,3)3}
            a_2a_3-a_1'b_1b_2 \equiv 0 \mod p.
        \end{equation}
    \end{enumerate}
    Thus, to summarize, we arrive at the following conditions
    \begin{enumerate}
    \item $a_1 = a_1'p$ where $0 \leq a_1' \leq p-1$,
    \item $(a_2^2 - a_2) - a_1'(b_1^2-b_1) \equiv 0 \mod p$,
    \item $a_3^2-b_2^2a_1' \equiv 0 \mod p$,
    \item $a_2a_3-a_1'b_1b_2 \equiv 0 \mod p.$
    \end{enumerate}
    We count matrices $A$ satisfying these conditions by dividing into cases.
    \par \textbf{Case 1 :} First consider the case when $a_3 \equiv \mod p$. From equation \eqref{(3,2,1,3)2} we deduce that $a_1'b_2 \equiv 0 \mod p$. We consider two subcases below.
    \begin{itemize}
        \item \textbf{Case 1A :} We consider the case when $b_2 = 0$. If $b_1 \in \{0,1\}$, then, there are $4p^3$ matrices. If $b_1 \notin \{0,1\}$, then, there are $p^3(p-2)$ matrices. Thus, there are $p^3(p+2)$ in this case in total.
        \item \textbf{Case 1B :} Next consider the case when $b_2 \neq 0$. From equation \eqref{(3,2,1,3)2} we find that $a_1' = 0$. From \eqref{(3,2,1,3)1}, we deduce $a_2^2-a_2 \equiv 0 \mod p$. Hence, there are $2p^3(p-1)$ matrices in this case.
    \end{itemize}
    \par \textbf{Case 2 :} In this case we count matrices such that $a_3 \not\equiv 0 \mod p$. From \eqref{(3,2,1,3)2}, we deduce that $b_2 \neq 0$. Consider following subcases.
    \begin{itemize}
        \item \textbf{Case 2A :} Suppose $b_1 = 0$, then from equation \eqref{(3,2,1,3)3} we get that $a_2 \equiv 0 \mod p$. Hence, there are $p^2(p-1)^2$ matrices in this case.
        \item \textbf{Case 2B :} Assume that $b_1 = 1$, then from equation \eqref{(3,2,1,3)1} we find that $a_2(a_2-1) \equiv 0 \mod p$. However, from \eqref{(3,2,1,3)2} we cannot have $a_2 \equiv 0 \mod p$, therefore $a_2 \equiv 1 \mod p$. Hence, we find that $a_3 \equiv b_2 \mod p$. As a result, there are $p^2(p-1)$ matrices in this case.
        \item \textbf{Case 2C :} In the last case we consider the case when $b_1 \notin \{0,1\}$ then from equations \eqref{(3,2,1,3)1}, \eqref{(3,2,1,3)2}, \eqref{(3,2,1,3)3} we get that $a_2 \not\equiv 0 \mod p$. Now, we get that from the three equations \eqref{(3,2,1,3)1}, \eqref{(3,2,1,3)2}, \eqref{(3,2,1,3)3} that,
        \begin{equation*}
            a_1' \equiv \frac{a_2(a_2-1)}{b_1(b_1-1)} \mod p,
        \end{equation*}
        \begin{equation*}
              a_1' \equiv \frac{a_3^2}{b_2^2} \mod p, 
        \end{equation*}
        \begin{equation*}
            a_1' \equiv \frac{a_2a_3}{b_1b_2} \mod p.
        \end{equation*}
        From these we find that $a_2 \equiv b_1 \mod p$ and $a_3 \equiv b_2 \mod p$. Hence, there are $(p-2)(p-1)p^2$ matrices in this case.
    \end{itemize}
    Adding up the values in each case, we prove the assertion of the lemma.
\end{proof}

\begin{lemma}\label{prev lemma}
    We have that $g_{\alpha}(p) = 2p^4$ where $\alpha = (3,1,3,2).$
\end{lemma}
\begin{proof}
    Let $A$ be any integer matrix in Hermite normal form of the type
    $$
    \begin{pmatrix}
        p^3 & a_1p & a_2p & a_3p & 1\\
         & p & 0 & 0 & 1\\
          & & p^3 & c_1p & 1\\
          &&& p^2 & 1\\
          &&&& 1\\
    \end{pmatrix}
    $$
    We derive the conditions so that $A$ is a subring matrix and use them to evaluate $g_{(3,1,3,2)}(p)$.
    \begin{enumerate}
        \item It is clear that $v_2^2 \in \Col(A)$ if and only if $v_2^2-pv_2 \in \Col(A)$, the second condition translates to
        \begin{equation}\label{(3,1,3,2)1}
            a_1^2-a_1 \equiv 0 \mod p.
        \end{equation}
        \item Similarly, $v_3^2 \in \Col(A)$ if and only if $v_3^2-p^3v_3 \in \Col(A)$, and therefore, $a_2 \equiv 0 \mod p$. In other words, $a_2 = a_2'p$, where $0 \leq a_2' \leq p-1$.
        \item Arguing as in previous two cases, we see that $v_4^2 \in \Col(A)$ if and only if $c_1 = c_1'p$, where $0 \leq c_1' \leq p-1$ and $a_3 = a_3'p$ with $0 \leq a_3' \leq p-1$.
        \item It is easy to see that if the above conditions are satisfied, then $v_2\circ v_3,v_2\circ v_4$ and $v_3\circ v_4$ are in $\Col(A)$.
    \end{enumerate}
    To summarize, the conditions are as follows
    \begin{enumerate}
        \item $a_1^2-a_1 \equiv 0 \mod p$,
        \item $a_2 = a_2'p$ where $0 \leq a_2' \leq p-1$,
        \item $c_1 = c_1'p$ where $0 \leq c_1' \leq p-1$, \item $a_3 = a_3'p$ with $0 \leq a_3' \leq p-1$.
    \end{enumerate}
    From the above, it is clear that $g_{(3,1,3,2)}(p) = 2p^4$.
\end{proof}
\begin{lemma}
    We have that
    \begin{equation}\label{(3,3,1,2)0}
        g_{(3,3,1,2)}(p)=p^5+4p^4-p^3.
    \end{equation}
\end{lemma}
\begin{proof}
Let $A$ be an integer matrix in Hermite normal form
$$
\begin{pmatrix}
p^3&a_1p&a_2p&a_3p&1\\
&p^3&b_1p&b_2p&1\\
&&p&0&1\\
&&&p^2&1\\
&&&&1
\end{pmatrix}.
$$
The condition $v_2^2\in\Col(A)$ gives $a_1=pa_1'$ with $0\leq a_1'\leq p-1$. The condition $v_3^2\in\Col(A)$ gives
\[b_1^2-b_1\equiv0\mod p\quad\text{and}\quad p(a_2^2-a_2)-a_1'(b_1^2-b_1)\equiv0\mod p^2.\]
Also, $v_4^2\in\Col(A)$ gives
\[b_2=pb_2'\quad\text{and}\quad a_3=pa_3',\]
where $0\leq b_2',a_3'\leq p-1$. Finally,
\[v_3\circ v_4\in\Col(A)\quad\text{if and only if}\quad a_1'b_1b_2'\equiv0\mod p.\]
These conditions are also sufficient.

Write $b_1=\epsilon+pt$ with $\epsilon\in\{0,1\}$ and $t\in\F_p$. Reducing the second congruence after division by $p$ gives
\[a_2(a_2-1)\equiv (2\epsilon-1)a_1't\mod p.\]
The residue of $a_2$ has one free lift modulo $p^2$, while $a_3'$ is free. If $\epsilon=0$, the condition involving $b_2'$ is automatic. The number of choices of $(a_1',t,a_2\bmod p)$ is $p^2+p$, while $b_2'$ is free. Including the free lift of $a_2$ and the free choice of $a_3'$, this case contributes
\[p^3(p^2+p)=p^5+p^4.\]
If $\epsilon=1$, the additional condition is $a_1'b_2'=0$. When $a_1'=0$ there are $2p^2$ residual choices for $(t,b_2',a_2\bmod p)$, while for $a_1'\neq0$ there are $p(p-1)$ residual choices. Including the two free lifts gives
\[p^2(3p^2-p)=3p^4-p^3.\]
Adding the two cases gives
\[g_{(3,3,1,2)}(p)=p^5+4p^4-p^3.\]
\end{proof}
\begin{lemma}
    The following equality holds
    \begin{equation}
        g_{(3,2,3,1)}(p)=2p^5+8p^4-8p^3+4p^2.
    \end{equation}
\end{lemma}
\begin{proof}
Let $A$ be an integer matrix in Hermite normal form of the type
$$
\begin{pmatrix}
p^3&a_1p&a_2p&a_3p&1\\
&p^2&b_1p&b_2p&1\\
&&p^3&c_1p&1\\
&&&p&1\\
&&&&1
\end{pmatrix}.
$$
As in the proof of Lemma \ref{(3,2,1,3)}, the condition $v_2^2\in\Col(A)$ gives $a_1=pa_1'$ with $0\leq a_1'\leq p-1$. The remaining square conditions imply
\[a_2^2-a_1'b_1^2\equiv0\mod p,\quad c_1(c_1-1)\equiv0\mod p,\quad\text{and}\quad b_1c_1(c_1-1)\equiv0\mod p^2.\]
The first-coordinate condition coming from $v_4^2$ and the mixed-product condition are
\begin{align*}
p^2(a_3^2-a_3)-a_2(c_1^2-c_1)-p^2a_1'(b_2^2-b_2)&\equiv0\mod p^3,\\
a_2(a_3-c_1)-a_1'b_1(b_2-c_1)&\equiv0\mod p.
\end{align*}

If $c_1=0$, the resulting system is the one occurring for the composition $(3,2,1,3)$, and it contributes
\[5p^4-4p^3+2p^2.\]
The same count holds when $c_1=1$. It remains to count the lifts for which $c_1\equiv0$ or $1\mod p$ but $c_1\notin\{0,1\}$ modulo $p^2$. In this case the condition $b_1c_1(c_1-1)\equiv0\mod p^2$ forces $b_1=0$, and then $a_2\equiv0\mod p$. Writing
\[c_1=\epsilon+pk\quad\text{where}\quad \epsilon\in\{0,1\}\quad\text{and}\quad k\in\F_p^\times,\]
the first-coordinate congruence determines the remaining lift of $a_2$ uniquely after $a_1',b_2$ and the residue of $a_3$ are chosen. The lift of $a_3$ is free. Thus every such $c_1$ contributes $p^4$ matrices. There are $2(p-1)$ choices for $c_1$, so these cases contribute
\[2p^4(p-1).\]
Therefore,
\begin{align*}
g_{(3,2,3,1)}(p)
&=2(5p^4-4p^3+2p^2)+2p^4(p-1)\\
&=2p^5+8p^4-8p^3+4p^2.
\end{align*}
\end{proof}
\begin{lemma}
    The following equality holds
    \[\begin{split}
        &g_{(3,1,2,3)}(p) = p^4 + 2p^3 - p^2,\\
        &g_{(3,3,2,1)}(p) = 6p^5-3p^4+4p^3,\\
        &g_{(3,2,2,2)}(p) = 2p^5-p^4.
    \end{split}\]
\end{lemma}
\begin{proof}
The proof is very similar to previous results, and we omit it.
\end{proof}

\subsection{Compositions of length $5$} We consider the composition of length $5$ that begin with $3$ we need to evaluate $g_{\alpha}(p)$ for $10$ compositions of this form. We list all compositions of this form in the table below\\
\begin{table}[H]
\centering
\begin{tabular}{c c } 
\hline
$(3,3,1,1,1)$ & $(3,2,1,1,2)$\\ 
$(3,1,3,1,1)$ & $(3,1,2,1,2)$\\
$(3,1,1,3,1)$ & $(3,1,2,2,1)$\\
$(3,1,1,1,3)$ & $(3,1,1,2,2)$\\
$(3,2,1,2,1)$ & $(3,2,2,1,1)$.\\
\hline
\end{tabular}
\caption{Compositions of length $5$ that begin with $3$.}
\end{table}
For some of these compositions, we can use results from \cite{isham2022number} to deduce value of $g_{\alpha}(p)$. The next result follows from the results in \emph{loc. cit.}
\begin{lemma}
We have the following values
\begin{table}[H]
\centering
 \begin{tabular}{c c } 
 \hline
 $g_{(3,3,1,1,1)}(p)$ & $16 p^6 + 12 p^5 - 20 p^4 + 8 p^3$\\ 
 $g_{(3,1,3,1,1)}(p)$ & $18 p^5 - 6 p^4$\\
 $g_{(3,1,1,3,1)}(p)$ & $2 p^5 + 10 p^4 - 4 p^3$\\
 $g_{(3,1,1,1,3)}(p)$ & $4p^4$.\\
\hline
\end{tabular}
\end{table}
\end{lemma}
\begin{proof}
    The above results are a direct consequence of \cite[Lemma 4.6]{isham2022number}.
\end{proof}
Next we evaluate the value of $g_{\alpha}(p)$ for the remaining compositions.
\begin{lemma}
The following table of relations holds
\begin{table}[H]
\centering
 \begin{tabular}{c c } 
 \hline
 $g_{(3,1,2,2,1)}(p)$ & $p^6+11p^5-6p^4$\\ 
 $g_{(3,1,2,1,2)}(p)$ & $6p^5-2p^4$\\
 $g_{(3,1,1,2,2)}(p)$ & $p^5+4p^4-2p^3$\\
 $g_{(3,2,1,1,2)}(p)$ & $4p^6+4p^5-9p^4+8p^3$\\
 $g_{(3,2,1,2,1)}(p)$ & $11p^6+p^5-15p^4+14p^3$.\\
\hline
\end{tabular}
\end{table}
\end{lemma}
\begin{proof}
We prove the assertion for the composition $(3,1,1,2,2)$; the other entries in the table are obtained in the same way. Let
\[A=\begin{pmatrix}
p^3&a_1p&a_2p&a_3p&a_4p&1\\
&p&0&0&0&1\\
&&p&0&0&1\\
&&&p^2&b_1p&1\\
&&&&p^2&1\\
&&&&&1
\end{pmatrix}.\]
The conditions $v_2^2,v_3^2\in\Col(A)$ give
\[a_1(a_1-1)\equiv0\mod p\quad\text{and}\quad a_2(a_2-1)\equiv0\mod p.\]
The condition $v_4^2\in\Col(A)$ gives $a_3=pa_3'$ with $0\leq a_3'\leq p-1$. For the fifth column, before substituting $a_3=pa_3'$, the first-coordinate condition is
\[p^2a_4^2-p^3a_4-p^2a_3b_1^2\equiv0\mod p^3,\]
and hence it becomes
\begin{equation}\label{(3,1,1,2,2) 0}
a_4^2-a_3'b_1^2\equiv0\mod p.
\end{equation}
The mixed products impose
\[a_1a_2\equiv0\mod p,\quad a_1a_4\equiv0\mod p,\quad\text{and}\quad a_2a_4\equiv0\mod p.\]
These conditions are also sufficient.

Let $\epsilon_i$ be the residue of $a_i$ modulo $p$ for $i=1,2$. Then
\[\epsilon_i\in\{0,1\}\quad\text{and}\quad \epsilon_1\epsilon_2=0.\]
There are $p^3$ free lifts coming from $a_1,a_2$ and $a_4$. If $\epsilon_1=\epsilon_2=0$, equation \eqref{(3,1,1,2,2) 0} has $p^2$ residual solutions in $(a_3',b_1,a_4\bmod p)$. If exactly one of $\epsilon_1,\epsilon_2$ is $1$, then $a_4\equiv0\mod p$, and there are $2p-1$ residual solutions in $(a_3',b_1)$. There are two choices for the nonzero $\epsilon_i$. Consequently,
\begin{align*}
g_{(3,1,1,2,2)}(p)
&=p^3\bigl(p^2+2(2p-1)\bigr)\\
&=p^5+4p^4-2p^3.
\end{align*}
\end{proof}

\par There is one composition $\alpha$ which requires a different argument, namely $\alpha=(3,2,2,1,1)$. We first prove the finite-field point count which occurs in this computation.

\begin{proposition}\label{prop:exceptional-point-count}
Let $N_p$ denote the number of solutions in $\F_p^8$ to
\[\begin{split}
&(x_3^2-x_3)-x_2(x_7^2-x_7)-x_1(x_5^2-x_5)=0,\\
&(x_4^2-x_4)-x_2(x_8^2-x_8)-x_1(x_6^2-x_6)=0,\\
&x_3x_4-x_2x_7x_8-x_1x_5x_6=0.
\end{split}\]
Then
\[N_p=p^5+12p^4-20p^3+30p^2-10p.\]
\end{proposition}
\begin{proof}
Set $h(x):=x(x-1)$ and
\[H:=\{(r_1,r_2,r_3)\in\F_p^3:r_1+r_2+r_3=1\}\quad\text{and}\quad Y(r_1,r_2,r_3):=(r_1r_2,r_1r_3,r_2r_3).\]
If $r=(x,y,1-x-y)$, then
\[(h(x),h(y),xy)=(-Y_1(r)-Y_2(r),-Y_1(r)-Y_3(r),Y_1(r)).\]
The linear transformation on the right is invertible. Therefore, after putting
\[r=(x_3,x_4,1-x_3-x_4),\quad s=(x_7,x_8,1-x_7-x_8),\quad\text{and}\quad t=(x_5,x_6,1-x_5-x_6),\]
the system in the statement is equivalent to
\begin{equation}\label{eq:Y-incidence}
Y(r)=x_2Y(s)+x_1Y(t).
\end{equation}
We count the solutions to \eqref{eq:Y-incidence} according to the dimension of the span of $Y(s)$ and $Y(t)$.

There are exactly three points of $H$ for which $Y(r)=0$, namely the three coordinate vertices; set $z:=3$. For a projective direction $d\in\mathbb{P}^2(\F_p)$, let $m_d$ be the number of $r\in H$ such that $Y(r)$ is nonzero and has direction $d$. The quadratic map
\[[r_1:r_2:r_3]\longmapsto[r_1r_2:r_1r_3:r_2r_3]\]
is the standard quadratic Cremona transformation. Hence $m_d=1$ when all three coordinates of $d=[u:v:w]$ are nonzero and
\[uv+uw+vw\neq0,\]
whereas the three coordinate directions have multiplicity $p-2$. All remaining directions have multiplicity zero. Thus, the number of directions of multiplicity $1$ is
\[p^2-3p+3.\]

If $Y(s)=Y(t)=0$, the contribution is
\begin{equation}\label{eq:C0}
C_0=z^3p^2=27p^2.
\end{equation}
Suppose next that $Y(s)$ and $Y(t)$ span a one-dimensional space. For a direction of multiplicity $m$, the number of ordered pairs $(s,t)$ on this direction, not both mapping to zero, is $(z+m)^2-z^2$. The map $(x_2,x_1)\mapsto x_2Y(s)+x_1Y(t)$ has fibers of size $p$, and the sum of the numbers of preimages under $Y$ along this line is $z+m$. It follows that
\begin{equation}\label{eq:C1}
\begin{split}
C_1={}&p(p^2-3p+3)((z+1)^2-z^2)(z+1)\\
&\quad +3p((z+p-2)^2-z^2)(z+p-2)\\
={}&3p^4+37p^3-102p^2+60p.
\end{split}
\end{equation}

It remains to treat the case when $Y(s)$ and $Y(t)$ are linearly independent. Let
\[\mathcal C:\ uv+uw+vw=0\]
be the corresponding nonsingular conic in $\mathbb P^2$. For a projective line $L$, set
\[M_L:=\sum_{d\in L}m_d\quad\text{and}\quad S_L:=\sum_{d\in L}m_d^2.\]
The contribution from ordered pairs whose images span $L$ is
\[(M_L^2-S_L)(z+M_L).\]
For $p\geq5$, the lines in $\mathbb P^2(\F_p)$ fall into the following six classes:
\begin{center}
\begin{tabular}{c|c|c|c}
class of $L$ & number & $M_L$ & $S_L$\\ \hline
coordinate side & $3$ & $2(p-2)$ & $2(p-2)^2$\\
tangent at a coordinate vertex & $3$ & $2p-3$ & $(p-1)+(p-2)^2$\\
other secant through a coordinate vertex & $3(p-2)$ & $2p-4$ & $(p-2)+(p-2)^2$\\
tangent at a noncoordinate point & $p-2$ & $p-3$ & $p-3$\\
secant through no coordinate vertex & $(p-2)(p-3)/2$ & $p-4$ & $p-4$\\
external line & $p(p-1)/2$ & $p-2$ & $p-2$
\end{tabular}
\end{center}
Summing over these classes gives
\begin{equation}\label{eq:C2}
C_2=p(p-2)(p^3+11p^2-35p+35)=p^5+9p^4-57p^3+105p^2-70p.
\end{equation}
For $p=2$ and $p=3$, a direct count gives $N_2=164$ and $N_3=915$, which agree with the formula in the statement. For $p\geq5$, adding \eqref{eq:C0}, \eqref{eq:C1} and \eqref{eq:C2} gives
\[N_p=p^5+12p^4-20p^3+30p^2-10p,\]
as required.
\end{proof}

\begin{theorem}\label{thm:exceptional-composition}
We have
\[g_{(3,2,2,1,1)}(p)=p^7+24p^6-29p^5+21p^4-4p^3.\]
\end{theorem}
\begin{proof}
Let $A$ be an integer subring matrix of the form
\[A=\begin{pmatrix}
p^3&a_1p&a_2p&a_3p&a_4p&1\\
&p^2&b_1p&b_2p&b_3p&1\\
&&p^2&c_1p&c_2p&1\\
&&&p&0&1\\
&&&&p&1\\
&&&&&1
\end{pmatrix}.\]
The condition $v_2^2\in\op{Col}(A)$ first gives $a_1=pa_1'$ with $0\leq a_1'\leq p-1$. We separate the count according as $b_1=0$ or $b_1\neq0$.

Suppose first that $b_1=0$. The remaining conditions force $a_2$ to be divisible by $p$. After removing the two free lifts, the residual equations are exactly the three equations in Proposition \ref{prop:exceptional-point-count}. Consequently, this case contributes
\[p^2N_p.\]

Suppose next that $b_1\neq0$. Writing $t\equiv a_2/b_1\mod p$, the square condition gives $a_1'=t^2$. The remaining conditions imply
\[c_1,c_2\in\{0,1\}\quad\text{and}\quad c_1c_2=0.\]
Writing $A,D$ for the residues of $a_3,a_4$ and $B,E$ for those of $b_2,b_3$, we obtain
\[h(A)=t^2h(B),\quad h(D)=t^2h(E),\quad\text{and}\quad AD=t^2BE,\]
together with
\[t(A-c_1-t(B-c_1))=0\quad\text{and}\quad t(D-c_2-t(E-c_2))=0.\]
If $t=0$, there are $9p^2$ residual solutions. If $t=1$, there are $3p^2$ residual solutions, while for each $t\notin\{0,1\}$ there are $3$ residual solutions. Hence the number of residual solutions is
\[9p^2+3p^2+3(p-2)=3(4p^2+p-2).\]
There are $p-1$ choices for $b_1$ and $p^3$ free lifts, so the contribution from $b_1\neq0$ is
\[3p^3(p-1)(4p^2+p-2).\]
Using Proposition \ref{prop:exceptional-point-count}, we conclude that
\begin{align*}
g_{(3,2,2,1,1)}(p)
&=p^2N_p+3p^3(p-1)(4p^2+p-2)\\
&=p^7+24p^6-29p^5+21p^4-4p^3.
\end{align*}
\end{proof}

\section{Calculating the values of $g_\alpha(p)$ for compositions beginning with $4,5$ or $6$}\label{s 6}
\subsection{Values of $g_\alpha(p)$ for compositions beginning with $4$}
\subsubsection{Compositions of length $4$}We will evaluate $g_{\alpha}(p)$ for $6$ compositions of length $4$ that begin with $4$. They are listed below\\
\begin{table}[h!]
\centering
\begin{tabular}{c c } 
\hline
$(4,2,2,1)$ & $(4,2,1,2)$\\ 
$(4,1,2,2)$ & $(4,3,1,1)$\\
$(4,1,3,1)$ & $(4,1,1,3)$.\\
\hline
\end{tabular}
\caption{Compositions of length $4$ that begin with $4$.}
\end{table}\\
\begin{lemma}
We have the following values,
\begin{table}[H]
\centering
 \begin{tabular}{c c } 
 \hline
 $g_{(4,2,2,1)}(p)$ & $13p^5-12p^4-2p^3+2p^2$\\ 
 $g_{(4,2,1,2)}(p)$ & $3p^5+4p^4-7p^3+p^2$\\
 $g_{(4,1,2,2)}(p)$ & $p^5 + 2 p^4 - p^2$\\
 $g_{(4,3,1,1)}(p)$ & $3p^5+15p^4-9p^3$\\
 $g_{(4,1,3,1)}(p)$ & $6p^4+2p^3-2p^2$\\
 $g_{(4,1,1,3)}(p)$ & $p^4+2p^3$.\\
\hline
\end{tabular}
\end{table}
\end{lemma}
\begin{proof}
The proof is omitted.
\end{proof}

\subsubsection{Compositions of length $5$}We will evaluate $g_{\alpha}(p)$ for $4$ compositions of length $5$ that begin with $4$. We list them below\\
\begin{table}[H]
\centering
\begin{tabular}{c c } 
\hline
$(4,2,1,1,1)$ & $(4,1,2,1,1)$\\ 
$(4,1,1,2,1)$ & $(4,1,1,1,2)$.\\
\hline
\end{tabular}
\caption{Compositions of length $5$ that begin with $4$.}
\end{table}
\begin{lemma}
We have the following values,
\begin{table}[H]
\centering
 \begin{tabular}{c c } 
 \hline
 $g_{(4,2,1,1,1)}(p)$ & $5p^6+22p^5-42p^4+16p^3$\\ 
 $g_{(4,1,2,1,1)}(p)$ & $17p^5-18p^4+3p^3$\\
 $g_{(4,1,1,2,1)}(p)$ & $5p^5+2p^4-4p^3$\\
 $g_{(4,1,1,1,2)}(p)$ & $p^5+3p^4$.\\
\hline
\end{tabular}
\end{table}
\end{lemma}
\begin{proof}
We first give the count for $(4,1,1,2,1)$, since this is the case in which the small-prime checks at $p=2$ and $p=3$ do not detect the error in the original expression. Let
\[A=\begin{pmatrix}
p^4&a_1p&a_2p&a_3p&a_4p&1\\
&p&0&0&0&1\\
&&p&0&0&1\\
&&&p^2&b_1p&1\\
&&&&p&1\\
&&&&&1
\end{pmatrix}.\]
The square conditions for the second and third columns give
\[a_i(a_i-1)\equiv0\mod p^2\quad\text{for}\quad i=1,2,\]
while the fourth-column condition gives $a_3=pc$ with $0\leq c<p^2$. The remaining conditions are
\begin{align*}
a_4(a_4-1)-c b_1(b_1-1)&\equiv0\mod p^2,\\
a_1a_2,\ a_1a_4,\ a_2a_4&\equiv0\mod p^2,\\
a_1c,\ a_2c,\ c(a_4-b_1)&\equiv0\mod p.
\end{align*}
Since $a_i(a_i-1)\equiv0\mod p^2$, we may write
\[a_i\equiv\epsilon_i\mod p^2\quad\text{where}\quad \epsilon_i\in\{0,1\}.\]
The condition $a_1a_2\equiv0\mod p^2$ shows that at most one of $\epsilon_1,\epsilon_2$ is equal to $1$.

If $\epsilon_1=\epsilon_2=0$, the free lifts of $a_1$ and $a_2$ contribute $p^2$. It remains to count triples $(a_4,c,b_1)$ satisfying
\[a_4(a_4-1)\equiv c b_1(b_1-1)\mod p^2\quad\text{and}\quad c(a_4-b_1)\equiv0\mod p.\]
If $c\equiv0\mod p$, this gives $2p^3$ triples. If $c\not\equiv0\mod p$, the cases $b_1\in\{0,1\}$ and $b_1\notin\{0,1\}$ give respectively $2p^2(p-1)$ and $p^2(p-2)$ triples. Thus this case contributes
\[p^2\bigl(2p^3+2p^2(p-1)+p^2(p-2)\bigr)=5p^5-4p^4.\]
If exactly one of $\epsilon_1,\epsilon_2$ is equal to $1$, there are two choices for its position. Then $a_4\equiv0\mod p^2$ and $c\equiv0\mod p$. For each position the count is $p^3(3p-2)$, and hence these cases contribute
\[2p^3(3p-2)=6p^4-4p^3.\]
Adding the contributions gives
\[g_{(4,1,1,2,1)}(p)=5p^5+2p^4-4p^3.\]

We now explicitly obtain the expression for $g_{(4,1,1,1,2)}(p)$. Let $A$ be a matrix in Hermite normal form of following type
$$
    \begin{pmatrix}
        p^4 & a_1p & a_2p & a_3p & a_4p & 1\\
         & p & 0 & 0 &  0 & 1\\
          & & p & 0 & 0 & 1\\
          &&& p & 0 & 1\\
          &&&& p^2& 1\\
          &&&&&1\\
    \end{pmatrix}.
    $$
    First, we determine the conditions on entries of the matrix which make $A$ a subring matrix. Note that for $i=2,3,4$, $v_i^2 \in \Col(A)$ if and only if $a_{i-1}(a_{i-1}-1) \equiv 0 \mod p^2$ and $v_5^2 \in \Col(A)$ if and only $a_4 \equiv 0 \mod p$. We also note that for $1<i<j$, $v_i\circ v_j \in \Col(A)$ if and only if $a_{i-1}a_{j-1} \equiv 0 \mod p^2$. First we suppose that for $i=2,3,4$, we have that $a_{i-1} \equiv 0 \mod p^2$. There are $p^5$ matrices of this form. Next, suppose exactly one of $a_1,a_2,a_3$ is congruent to $1$ modulo $p^2$ and the others are divisible by $p^2$. Then we get $3p^4$ matrices of this form. Therefore,
    \[g_{(4,1,1,1,2)}(p)=p^5+3p^4.\]
\end{proof}

\subsection{Values of $g_\alpha(p)$ for compositions beginning with $5$}
\subsubsection{Compositions of length $4$}We calculate $g_{\alpha}(p)$ for $3$ compositions of length $4$ that begin with $5$. The compositions are listed below\\
\begin{table}[H]
\centering
\begin{tabular}{c } 
\hline
$(5,2,1,1)$\\ 
$(5,1,2,1)$\\
$(5,1,1,2)$.\\
\hline
\end{tabular}
\caption{Compositions of length $4$ that begin with $5$.}
\end{table}
\begin{lemma}
The following table of relations holds
\begin{table}[H]
\centering
 \begin{tabular}{c c } 
 \hline
 $g_{(5,2,1,1)}(p)$ & $19p^4-21p^3+6p^2$\\
 $g_{(5,1,2,1)}(p)$ & $7 p^4 - p^3 - 2 p^2$\\
 $g_{(5,1,1,2)}(p)$ & $2p^4+2p^3$.\\
\hline
\end{tabular}
\end{table}
\end{lemma}

\begin{proof}
We will prove the result only for the composition $(5,1,2,1)$. The other cases are omitted since the arguments are similar to this case. Let $A$ be a matrix in Hermite normal form of the following type
\[
\begin{pmatrix}
        p^5 & a_1p & a_2p & a_3p & 1\\
         & p & 0 & 0 & 1\\
          & & p^2 & b_1p & 1\\
          &&& p & 1\\
          &&&& 1\\
\end{pmatrix}\]
    We obtain the conditions on the entries of $A$ below.
    \begin{enumerate}
        \item First, we must have that $v_2^2 \in \Col(A)$ and this is true if and only if $v_2^2-pv_2 \in \Col(A)$. The latter condition is equivalent to the condition $a_1(a_1-1) \equiv 0 \mod p$.
        \item We also want that $v_3^2 \in \Col(A)$ and we see that this is equivalent to the condition $v_3^2-p^2v_3 \in \Col(A)$, which is true if and only if $a_2(a_2-p) \equiv \mod p$. Therefore we get that $a_2 = a_2'p$ with $0 \leq a_2' \leq p^3-1$ and $a_2'(a_2'-1) \equiv 0 \mod p$.
        \item Using similar arguments we obtain that $v_4^2 \in \Col(A)$ if and only if $a_3(a_3-1)-b_1(b_1-1)a_2' \equiv 0 \mod p^3$. We also deduce that $v_2\circ v_3, v_2\circ v_4$ and $v_3\circ v_4$ are in $\Col(A)$ if and only if $a_1a_2 \equiv 0 \mod p$, $a_1a_3 \equiv 0 \mod p$ and $a_2(a_3-b_1) \equiv 0 \mod p$.
    \end{enumerate}
    The conditions we get are as follows
    \begin{itemize}
        \item $a_1(a_1-1) \equiv 0 \mod p$,
        \item $a_2 = a_2'p$ with $0 \leq a_2' \leq p^3-1$ and $a_2'(a_2'-1) \equiv 0 \mod p$,
        \item $a_3(a_3-1)-b_1(b_1-1)a_2' \equiv 0 \mod p^3$,
        \item $a_1a_2 \equiv 0 \mod p$, $a_1a_3 \equiv 0 \mod p$ and $a_2(a_3-b_1) \equiv 0 \mod p$.
    \end{itemize}
    Next, we count the number of matrices by dividing into two cases according as $a_1 \equiv 1 \mod p^3$ or $a_1 \equiv 0 \mod p^3$.
    \par \textbf{Case 1 :} First consider the case when $a_1 \equiv 1 \mod p^3$, as we must have $a_1a_2 \mod p^3$ we deduce that $a_2 = a_2'''p^3$ and $a_3 \equiv 0 \mod p$. Therefore, conditions reduce to $b_1(b_1-1)a_2''' \equiv 0 \mod p$. Counting we see that there are $3p^3-2p^2$ matrices in this case.
    \par \textbf{Case 2 :} Now, we consider that case when $a_1 \equiv 0 \mod p^3$. We further consider two separate cases depending upon whether $a_2'$ is divisible by $p$ or not.
    \begin{itemize}
        \item \textbf{Case 2A :} Consider the case when $a_2'$ is divisible by $p$. we have two possibilities, either $p^3 \nmid a_2$ or $p^3 \mid a_2$, if $p^3 \nmid a_2$ then we deduce that $a_3 \equiv b_1 \mod p$, therefore $b_1 = 0$ or $1$. We get $2p^3(p-1)$ matrices. if $p^3 \mid a_2$, again considering two cases depending upon whether $b_1 \in \{0,1\}$ or not and counting we see that we get $2p^4$ matrices.
        \item \textbf{Case 2B :} We consider the case when $a_2' \equiv 1 \mod p$, we deduce that $a_3 \equiv b_1 \mod p^3$. Conditions reduce to $b_1(b_1-1)(1-a_2') \equiv 0 \mod p^3$. Now, suppose $b_1 \in \{0,1\}$, then we get $2p^4$ matrices. Otherwise, $a_2' = 1$ and we get $p^2(p-2)$ matrices. 
    \end{itemize}
    It follows from the above observations that the number of matrices in each case we deduce that $g_{(5,1,2,1)}(p) = 7 p^4 - p^3 - 2 p^2$.
\end{proof}

\subsubsection{Compositions of length $5$}In this section we evaluate $g_{\alpha}(p)$ for the composition $\alpha = (5,1,1,1,1)$. This is the only composition of length $5$ which starts with $5$. We have the following result.
\begin{lemma}
We have that $g_{\alpha}(p) = 5p^4$, where $\alpha = (5,1,1,1,1)$.
\end{lemma}
\begin{proof}
The result follows directly from \cite[Lemma 3.5]{atanasov2021counting}.
\end{proof}

\subsection{Values of $g_\alpha(p)$ for compositions beginning with $6$}
We evaluate $g_{\alpha}(p)$ for the only composition $\alpha = (6,1,1,1)$ of length $4$ which begins with $6$.
\begin{lemma}
We have that $g_{(6,1,1,1)}(p) = 4p^3$.
\end{lemma}
\begin{proof}
The above result is an immediate consequence of \cite[Lemma 3.5]{atanasov2021counting}.
\end{proof}

\section{Main results}
\subsection{Proof of the main result}\label{s 7.1}In this section we prove Theorem \ref{main}. We note that for $n > 1$ and $e \geq n-1$ the following relation holds
\begin{equation}\label{eq7.1}
    g_n(p^e) = g_{n-1}(p^{e-1}) + \sum_{\alpha \in \mathcal{C}'_{n,e}}g_{\alpha}(p),
\end{equation}
where $\mathcal{C}'_{n,e}$ denotes the set of compositions in $\mathcal{C}_{n,e}$ whose first coordinate is greater than $1$.
\begin{proposition}
    We have that
    \[\begin{aligned}
    g_5(p^9)&=36p^5+126p^4-33p^3+12p^2+p+1,\\
    &\quad\text{and}\quad
    g_6(p^9)=p^7+112p^6+133p^5-76p^4+48p^3+2p^2+p+1.
    \end{aligned}\]
\end{proposition}
\begin{proof}
The result follows immediately from \eqref{eq7.1} and computations in previous sections.
\end{proof}
Liu computes the values of $g_n(p^e)$ for $n=3,4$ (cf. \cite[Proposition 6.1 and 6.2]{liu2007counting}). We record these values for $e=9$ and refer to \cite[p. 231]{atanasov2021counting} the values of $e$ which lie in the range $4\leq e\leq 8$. We find that
\begin{equation*}
    g_3(p^9) = p^3 + 4p^2 + 4p + 1,
\end{equation*}
and that
\begin{equation*}
    g_4(p^9) = 11p^4+30p^3+9p^2+p+1.
\end{equation*}

\begin{proposition}\label{mainprop}
For $n > 1$, the following relation holds
\vspace{-0.2 em}
\begin{eqnarray*}
    f_{n}(p^9) - f_{n-1}(p^9) & = & Q(n)(276480P_{10}(n)p^{10} + 276480P_{9}(n)p^{9} + 138240P_{8}(n)p^{8}\\
    & & + 34560P_{7}(n)p^{7} + 34560P_{6}(n)p^{6} + 11520P_{5}(n)p^{5}\\
    & & + 11520P_{4}(n)p^{4} + 1440P_{3}(n)p^{3} + 1440P_{2}(n)p^{2}\\
    & & + 240P_{1}(n)p + P_{0}(n)) + {n-1 \choose 5}g_{(3,2,2,1,1)}(p)\\
    & & + {n-1 \choose 4}\Delta_5(p)+{n-1 \choose 5}\Delta_6(p),
\end{eqnarray*}
where
\[\begin{aligned}
\Delta_5(p)&:=-11p^6+22p^5-11p^4-17p^3+11p^2-p-2,\\
&\quad\text{and}\quad
\Delta_6(p):=-4p^7+12p^6+34p^5-153p^4+163p^3-41p^2+10p.
\end{aligned}\]
Moreover,
\begin{eqnarray*}
    Q(n) &=& \frac{n-1}{1393459200},\\
    P_{10}(n) &=& (n - 2)(n - 3)(n - 4)(n - 5)(n - 6)(n - 7),\\
    P_9(n) &=& n(n - 2)(n - 3)(n - 4)(n - 5)(n - 6),\\
    P_8(n) &=& (n - 2)(n - 3)(n - 4)(n - 5)(n - 6)(9n^2 - 131n + 784),\\
    P_7(n) &=& (n - 2)(n - 3)(n - 4)(n - 5)(37n^3 - 761n^2 + 6482n - 18144),\\
    P_6(n) &=& (n - 2)(n - 3)(n - 4)(75n^5 - 2468n^4 + 36349n^3 - 279672n^2\\
    && + 1101372n - 1732080).\\
    P_5(n) &=& (n - 2)(n - 3)(n - 4)(33n^6 - 1220n^5 + 21757n^4 - 208732n^3\\
    && + 1053542n^2 - 2414076n + 1592640),\\
    P_4(n) &=& (n - 2)(n - 3)(42n^8 - 2102n^7 + 51106n^6 - 749135n^5\\
    && + 7010050n^4 - 41819711n^3 + 152633830n^2\\
    && - 307862664n + 260890560),\\
    P_3(n) &=& (n - 2)(273n^{10} - 18123n^9 + 571282n^8 - 10999886n^7\\
    && + 140978185n^6 - 1241293667n^5 + 7531038196n^4\\
    && - 30849468932n^3 + 81151302432n^2 - 123180801984n\\
    && + 81613163520),\\
    P_2(n) &=& (n - 2)(147n^{11} - 10843n^{10} + 377832n^9 - 8095642n^8\\
    && + 117343063n^7 - 1198590955n^6 + 8745674590n^5\\
    && - 45324367680n^4 + 162663439888n^3 - 383226514464n^2\\
    && + 531138427776n - 326776343040),\\
    P_1(n) &=& (n - 2)(315n^{12} - 25368n^{11} + 956025n^{10} - 22147762n^9\\
    && + 349611873n^8 - 3946086924n^7 + 32542876455n^6\\
    && - 196940866290n^5 + 865325849028n^4 - 2683887407672n^3\\
    && + 5560471449216n^2 - 6887555482752n + 3844971970560),\\
                 \end{eqnarray*}
    \begin{eqnarray*}
    P_0(n) &=& 135 n^{16} - 14400 n^{15} + 730980 n^{14} - 23334360 n^{13} + 522513250 n^{12}\\
    && - 8678453720 n^{11} + 110319301164 n^{10} - 1092302298312 n^9\\
    && + 8494128305343 n^8 - 51930963880392 n^7 + 248190356069720 n^6\\
    && - 915171074718208 n^5  + 2545385435375472 n^4\\
    && - 5146605000021888 n^3 + 7110255039457536 n^2\\
    && - 5973926003435520 n + 2289569122713600.
\end{eqnarray*}
\end{proposition}
\begin{proof}
\par Recall the recurrence relation \eqref{recursive formula}
\[f_n(p^9)=\sum_{i=0}^9\sum_{j=1}^n {n-1\choose j-1} f_{n-j}(p^{9-i})g_j(p^i).\] Note that since $i\leq 9$, we find that $j\leq 10$ in the above recurrence relation. We compute the values of $g_j(p^i)$ and then recursively, use the above to compute the value of $f_n(p^9)$. In greater detail, we use values from \cite[p. 231]{atanasov2021counting}, the Theorem \ref{theorem1.2}, and the values computed in previous sections, to obtain the above recurrence relation.
\end{proof}
Using the recurrence relation above and noting that \[f_n(p^9)=\sum_{k=2}^{n}\left(f_k(p^9)-f_{k-1}(p^9)\right)\]we obtain Theorem \ref{main}. We performed the symbolic computations in \begin{tt}SageMath\end{tt}. In addition, we independently checked the displayed formulae for $g_\alpha(p)$ by exact enumeration of the corresponding Hermite normal form matrices for $p=2$ and $p=3$. The corrected exceptional rows were also checked at $p=5$, and, whenever the size of the enumeration was reasonable, at $p=7$.
\subsection{Some bounds on $g_\alpha(p)$}\label{s 7.2}In this section we obtain bounds on $g_{\alpha}(p)$ for certain compositions $\alpha$. Let $u_1, \dots, u_n$ be the standard basis of $\Z^n$. The vector $u_i=(0,0,\dots, 0, 1, 0,\dots, 0)$ consists of $0$ in all entries except for $1$ in the $i$-th entry. Let $L$ be an irreducible subring of $\Z^n$ and let $\mathfrak{m}_L$ be the ideal in $L$ consisting of vectors all of whose coordinates are divisible by $p$. Setting $\rho_L:=\op{dim}_{\F_p} \mathfrak{m}_L/\mathfrak{m}_L^2$. Following \cite[p. 292, l.-1]{liu2007counting}, an irreducible subring is \emph{full} if $\rho_L=n-1$. Assume that the matrix for $L$ is in Hermite normal form with respect to the basis $u_1, \dots, u_n$. Let $\pi:\Z^n\rightarrow \Z^{n-1}$ be the projection onto the last $(n-1)$ coordinates $u_2, \dots, u_n$. Suppose that the basis is chosen so that with respect to the ordered basis $(u_1, \dots, u_n)$, the matrix $A$ associated with $L$ is in the Hermite normal form. Let $v_1, \dots, v_n$ be the columns of $A$, and let $\tau(L)$ be the lattice generated by $L$ and by $v_1':=\frac{1}{p} v_1$. It is easy to see that $\tau(L)$ is a subring (cf. \cite[p.291, l.-2]{liu2007counting}).
\begin{proposition}\label{liu prop5.6}
 Let $L$ be an irreducible subring of index $p^{e+n}$ such that $\pi(L)$ is full and $pu_1 \notin L$. Then the map $L \mapsto \frac{1}{p}\mathfrak{m}_{\tau(L)}$ is a $p^{n-2}$-to-one surjection onto subrings of index $p^e$.
\end{proposition}
\begin{proof}
The above result is \cite[Proposition 5.6]{liu2007counting}.
\end{proof}
Let $\beta := (\beta_1,\beta_2,\dots,\beta_{n-1})$ be a composition of length $n-1$ such that $\beta_1 > 1$, $\beta_i > 0$ for $i \in \{2,\dots,n-1\}$. By abuse of notation, we say that an irreducible subring $L$ has diagonal $\beta$ if the entries on the diagonal of the matrix $A_L$ are $p^{\beta_1}, \dots, p^{\beta_{n-1}},1$. We set $\mathcal{R_{\beta}}$ to denote the set of irreducible subrings with diagonal $\beta$. Let $\beta' := (\beta_1-2,\beta_2-1,\dots,\beta_{n-1}-1)$, let $\mathcal{S}_{\beta'}$ denote the set of subring matrices with diagonal $\beta'$. Let $\mathcal{R}'_\beta$ be the subset of $\mathcal{R}_\beta$ consisting of irreducible subrings with diagonal $\beta$ such that $\pi(L)$ is full and $pu_1\notin L$, with respect to the basis $u_1, \dots, u_n$ for which $A_L$ is in Hermite normal form. Given a finite set $S$, set $\# S$ to denote the cardinality of $S$. 
\begin{lemma}\label{7.2}
Let $\beta := (\beta_1,\beta_2,\dots,\beta_{n-1})$ be a composition of length $n-1$ such that the following conditions are satisfied
\begin{enumerate}
    \item $\beta_1 > 1$, 
    \item $\beta_i > 0$ for $i \in \{2,\dots,n-1\}$,
    \item $\sum_{i=1}^{n-1}\beta_i > n$.
\end{enumerate} Then, there is a $p^{n-2}$-to-one surjection from $\mathcal{R}_{\beta}'$ to $\mathcal{S}_{\beta'}$. In particular, we have that
\[\#\mathcal{R}_{\beta} \geq \#\mathcal{R}'_{\beta}\geq p^{n-2}\#\mathcal{S}_{\beta'}.\]
\end{lemma}
\begin{proof}
It follows from Proposition \ref{liu prop5.6} that the map $L \mapsto \frac{1}{p}\mathfrak{m}_{\tau(L)}$ defines a $p^{n-2}$ to one surjective map $\Phi:\mathcal{R}_\beta'\rightarrow \mathcal{S}_{\beta'}$. The inequality follows immediately from this. 
\end{proof}
The following result comes as a consequence of the above Lemma. 
\begin{proposition}\label{prop 7.3}
Let $\beta := (\beta_1,\beta_2,\dots,\beta_{n-1})$ be a composition of length $n-1$ such that $\beta_1 > 1$, $\beta_i > 0$ for $i \in \{2,\dots,n-1\}$ and $\sum_{i=1}^{n-1}\beta_i > n$. Setting $m_\beta := \min\{\lfloor \frac{\beta_1}{2}\rfloor,\beta_2,\dots,\beta_{n-1}\}$, we have the following lower bound for $g_\beta(p)$
\[g_{\beta}(p) \geq p^{(n-2)m_\beta}.\]
\end{proposition}
\begin{proof}
   The result follows upon repeatedly applying the Lemma \ref{7.2}.
\end{proof}
\begin{remark}
The lower bound in Proposition \ref{prop 7.3} is not expected to be sharp in general. For example, Theorem \ref{thm:exceptional-composition} gives a polynomial of degree $7$ for $g_{(3,2,2,1,1)}(p)$, whereas Proposition \ref{prop 7.3} only gives the lower bound $p^4$. More generally, after $r$ applications of Lemma \ref{7.2}, one retains the stronger intermediate estimate
\[g_\beta(p)\geq p^{r(n-2)}\#\mathcal S_{\beta^{(r)}},\quad\text{where}\quad \beta^{(r)}=(\beta_1-2r,\beta_2-r,\dots,\beta_{n-1}-r).\]
Thus, any additional lower bound for the terminal quantity $\#\mathcal S_{\beta^{(r)}}$ gives a corresponding improvement. Optimizing this terminal contribution over compositions is a natural possible refinement of the argument.
\end{remark}
In particular we see that if $\alpha := (k,\ell,\dots,\ell)$ of length $n-1$, then $g_{\alpha}(p) \geq p^{{\lfloor\frac{k}{2}\rfloor}(n-2)}$. We note that this is a special case of \cite[Corollary 4.6]{ishamlowerbound}. 
\begin{theorem}
    We have the following polynomial lower bound 
    \[g_n(p^e)\geq \sum_{j=1}^{\lfloor e/n\rfloor} p^{(n-2)j}\left( \binom{e-1-nj}{n-2}-\binom{e-1-n(j+1)}{n-2}\right).\]
\end{theorem}
\begin{proof}
Let $\mathcal{C}_{n,e}^{j}$ be the subset of $\mathcal{C}_{n,e}$ consisting of all compositions $\beta=(\beta_1, \dots, \beta_{n-1})$ such that $m_\beta=j$. It is easy to see that 
\[\#\mathcal{C}_{n,e}^j=\# \mathcal{C}_{n, e-nj}-\# \mathcal{C}_{n, e-n(j+1)}=\binom{e-1-nj}{n-2}-\binom{e-1-n(j+1)}{n-2}.\]

It follows from Proposition \ref{prop 7.3} that
\[g_n(p^e)\geq \sum_{j=1}^{\lfloor e/n\rfloor} p^{(n-2)j} \# \mathcal{C}_{n,e}^j,\]and the result follows.
\end{proof}

\par We now prove Theorem \ref{main2}.
\begin{proof}[Proof of Theorem \ref{main2}]
\par Let us assume for simplicity that $t=n-1$. The proof in generality is identical to this case. Given an $n\times n$ matrix $A$, let $A_{[i,j]}$ be the $i\times j$ matrix obtained upon deleting the last $(n-i)$ rows and $(n-j)$ columns. Thus, $A_{[i,j]}$ is the upper left $i\times j$ submatrix of $A$. Let $S$ denote the set of subring matrices of the form
$$A = \begin{pmatrix}
p^{\alpha_1} & \dots & \dots & \dots & a_1p & 1\\
 & p^{\alpha_2} & \dots & \dots & a_2p & 1\\
 &  & \ddots &  & \vdots & \vdots \\ 
 &  &   & & \vdots & \vdots\\
 &  &  &  & p^{k\alpha_{n-1}} & 1\\
 & &   &  &  & 1
\end{pmatrix}$$
that are in Hermite normal form. We note that the following conditions hold
\begin{enumerate}
    \item $0 \leq a_i \leq p^{\gamma-1}-1$,
    \item setting $\hat{A}:=A_{[n-2, n-2]}$, we note that \[(a_1^2p^2-a_1p^{k\alpha_{n-1}+1}, \dots, a_{n-2}^2p^2 - a_{n-2}p^{k\alpha_{n-1}+1})^t \in \op{Col}\hat{A}=\hat{A}(\Z^{n-2}).\]
\end{enumerate}
We note in passing that $g_{\alpha_{/k,t}}(p)= \#S$. We denote the set of all irreducible subring matrices of size $(n-1)$ and with diagonal $(p^{\alpha_1},p^{\alpha_2},\dots,p^{\alpha_{n-2}}, 1)$ by $\mathcal{A}$. 
\par We note that for each index $i$ in the range $1 \leq i \leq n-2$ we have that $|a_i^2p^2-a_ip^{k\alpha_i+1}| \leq |a_i^2|p^2 + |a_ip^{k\alpha_i+1}| \leq 2p^{(k+1)\gamma}$. 

\par In order to simplify notation, for $B\in \mathcal{A}$, set $B':=B_{[n-2, n-2]}$. We also set
\[\mathcal{C}:=[-2p^{(k+1)\gamma},2p^{(k+1)\gamma}]^{n-2}.\]
We then find that
\begin{equation}\label{ineq1}
    \#S \leq \sum_{B \in \mathcal{A}}\#\left(B'(\Z^{n-2})\cap \mathcal{C}\right).
\end{equation}
For $B \in \mathcal{A}$, the number of points in $v \in \Z^{n-2}$ such that $B'v \in \mathcal{C}$ is same as the number of points $w \in \mathcal{C}$ such that $\left(B'\right)^{-1}(w) \in \Z^{n-2}$. Therefore using inequality \ref{ineq1} we deduce that
\begin{equation}\label{ineq2}
    \#S \leq \sum_{B\in \mathcal{A}}\#(\left(B'\right)^{-1}(\mathcal{C})\cap \Z^{n-2}).
\end{equation}
For $w=(w_1, \dots, w_{n-2})\in \Z^{n-2}$, set $\lVert w \rVert$ to be the Euclidean norm $\sqrt{\sum_i w_i^2}$. For any $w \in \mathcal{C}$ we have that $\lVert (B')^{-1}w \rVert \leq \lVert (B')^{-1} \rVert\lVert w \rVert$. It is easy to see that $\lVert w \rVert \leq 2p^{(k+1)\gamma}\sqrt{n-2}$. On the other hand, since $\mathcal{A}$ is a finite set which is defined independent of $k$, we find that for any $B \in \mathcal{A}$ and $w \in \mathcal{C}$, $\lVert (B')^{-1}w \rVert \leq M\sqrt{n-2}p^{(k+1)\gamma}$ where $M \in \mathbb{R}_{>0}$ is a suitably large constant (not depending on $k$). Using the inequality \eqref{ineq2}, we deduce that
\begin{equation}\label{ineq3}
    \#S \leq \sum_{L \in \mathcal{A}}\#([-M\sqrt{n-2}p^{(k+1)\gamma},M\sqrt{n-2}p^{(k+1)\gamma}]^{n-2}\cap \Z^{n-2})
\end{equation}
Fixing an integer $N > M\sqrt{n-2}$, we deduce that
\begin{equation}\label{ineq4}
    \#S \leq \sum_{L \in \mathcal{A}}\#([-Np^{(k+1)\gamma},Np^{(k+1)\gamma}]^{n-2}\cap \Z^{n-2})
\end{equation}
Therefore, $\#S \leq \#(\mathcal{A})(2N)^{n-2}p^{(k+1)\gamma {(n-2)}}$. We have proved that
\[g_{\alpha_{/k,t}}(p) =\#S \leq \#(\mathcal{A})(2N)^{n-2}p^{\gamma{(n-2)}}p^{\gamma k{(n-2)}}\]
Therefore, $g_{\alpha_{/k,t}}(p) = O(p^{\gamma k(n-2)})\text{ as }k \to \infty.$
\end{proof}
\bibliographystyle{alpha}
\bibliography{references}

@article{grunewald1988subgroups,
  author  = {Grunewald, Fritz J. and Segal, Daniel and Smith, Geoff C.},
  title   = {Subgroups of finite index in nilpotent groups},
  journal = {Inventiones Mathematicae},
  volume  = {93},
  year    = {1988},
  pages   = {185--223},
  doi     = {10.1007/BF01393692}
}

@book{du2008zeta,
  author    = {du Sautoy, Marcus and Woodward, Luke},
  title     = {Zeta Functions of Groups and Rings},
  series    = {Lecture Notes in Mathematics},
  volume    = {1925},
  publisher = {Springer-Verlag},
  address   = {Berlin},
  year      = {2008}
}

@book{lubotzky2003subgroup,
  author    = {Lubotzky, Alexander and Segal, Dan},
  title     = {Subgroup Growth},
  series    = {Progress in Mathematics},
  volume    = {212},
  publisher = {Birkh\"auser Verlag},
  address   = {Basel},
  year      = {2003}
}

@phdthesis{brakenhoff2009counting,
  author = {Brakenhoff, Jos F.},
  title  = {Counting Problems for Number Rings},
  school = {Leiden University},
  year   = {2009}
}

@article{liu2007counting,
  author  = {Liu, Ricky Ini},
  title   = {Counting subrings of {$\mathbb Z^n$} of index {$k$}},
  journal = {Journal of Combinatorial Theory, Series A},
  volume  = {114},
  number  = {2},
  year    = {2007},
  pages   = {278--299}
}

@article{atanasov2021counting,
  author  = {Atanasov, Stanislav and Kaplan, Nathan and Krakoff, Benjamin and Menzel, Julia H.},
  title   = {Counting finite index subrings of {$\mathbb Z^n$}},
  journal = {Acta Arithmetica},
  volume  = {197},
  number  = {3},
  year    = {2021},
  pages   = {221--246},
  doi     = {10.4064/aa180201-29-7}
}

@article{kaplan2015distribution,
  author  = {Kaplan, Nathan and Marcinek, Jake and Takloo-Bighash, Ramin},
  title   = {Distribution of orders in number fields},
  journal = {Research in the Mathematical Sciences},
  volume  = {2},
  year    = {2015},
  pages   = {Article 6},
  doi     = {10.1186/s40687-015-0027-8}
}

@article{datskovsky1988density,
  author  = {Datskovsky, Boris and Wright, David J.},
  title   = {Density of discriminants of cubic extensions},
  journal = {Journal f\"ur die reine und angewandte Mathematik},
  volume  = {386},
  year    = {1988},
  pages   = {116--138}
}

@book{nakagawa1996orders,
  author    = {Nakagawa, Jin},
  title     = {Orders of a quartic field},
  series    = {Memoirs of the American Mathematical Society},
  volume    = {122},
  publisher = {American Mathematical Society},
  year      = {1996},
  note      = {No. 583, viii+75 pp.}
}

@article{ishamlowerbound,
  author  = {Isham, Kelly},
  title   = {Lower bounds for the number of subrings in {$\mathbb Z^n$}},
  journal = {Journal of Number Theory},
  volume  = {234},
  year    = {2022},
  pages   = {363--390},
  doi     = {10.1016/j.jnt.2021.06.026}
}

@article{isham2022number,
  author  = {Isham, Kelly},
  title   = {Is the number of subrings of index {$p^e$} in {$\mathbb Z^n$} polynomial in {$p$}?},
  journal = {Communications in Algebra},
  volume  = {51},
  number  = {12},
  year    = {2023},
  doi     = {10.1080/00927872.2023.2225600}
}

@article{hausel2008mixed,
  author  = {Hausel, Tam\'as and Rodriguez-Villegas, Fernando},
  title   = {Mixed Hodge polynomials of character varieties},
  journal = {Inventiones Mathematicae},
  volume  = {174},
  number  = {3},
  year    = {2008},
  pages   = {555--624},
  doi     = {10.1007/s00222-008-0142-x}
}
\end{document}